\documentclass[12pt,x11names]{article}
\include{drlmacros}
\usepackage{hhline}
\usepackage{comment}
\numberwithin{equation}{section}
\newcommand{\ie}{{\em i.e.}, }

\newcommand{\nn}{\mathbb{N}} 
\newcommand{\rr}{\mathbb{R}} 

\newcommand{\bu}{{\bf u}}

\newcommand{\by}{{\bf y}}
\newcommand{\bz}{{\bf z}}

\title{Optimization on Spheres: Models and Proximal Algorithms with Computational Performance Comparisons}
\author{D.\ Russell Luke\thanks{Institut f\"ur Numerische und Angewandte Mathematik, Universit\"at G\"ottingen,\ Lotzestr.~16--18, 37083 G\"ottingen, Germany. E-mail: \texttt{r.luke@math.uni-goettingen.de}. Supported in part by the Deutsche Forschungsgemeinschaft grant SFB755 and by GIF Grant G-1253-304.6/2014.} \and Shoham Sabach\thanks{Faculty of Industrial Engineering and Management, Technion---Israel Institute of Technology, Haifa 3200003, Israel. E-mail: ssabach@ie.technion.ac.il. Partially supported by GIF Grant G-1253-304.6/2014.} \and Marc Teboulle\thanks{School of Mathematical Sciences, Tel-Aviv University, Ramat-Aviv 69978, Israel. E-mail: teboulle@post.tau.ac.il. Partially supported by the Israel Science Foundation, ISF Grant 998/12 and by GIF Grant G-1253-304.6/2014.}}
\date{9.06.2018}

\begin{document}
\maketitle

\begin{abstract}
We present a unified treatment of the abstract problem of finding the best approximation 
between a cone and spheres in the image of affine transformations.  Prominent instances 
of this problem are phase retrieval and source localization.  The common geometry binding 
these problems permits a generic application of algorithmic ideas and abstract convergence 
results for nonconvex optimization.  We organize  variational models for this problem into 
three different classes and derive 
the main algorithmic approaches within these classes (13 in all).  We identify the central ideas 
underlying these methods and provide thorough numerical 
benchmarks comparing their performance on synthetic and laboratory data. The software 
and data of our experiments are all publicly accessible. We also introduce one new algorithm, 
a cyclic relaxed Douglas-Rachford algorithm, which outperforms all other algorithms by every 
measure: speed, stability and accuracy. The analysis of this algorithm remains open.
\end{abstract}




\section{Introduction}
Nonconvex optimization is maturing into a major focus within continuous optimization. The numerical intuition that one brings to this area from convex optimization can be misleading. The purpose of this work is to provide thorough numerical benchmarks on a prevalent type of nonconvex problem which will establish guideposts for the good, the bad and the ugly of numerical methods. The problem we consider is the following:
    \begin{center}
    		\fcolorbox{black}{Ivory2}{\parbox{15cm}{\vspace{-0.0in}
        		\textbf{The Cone and Sphere Problem:}\\ Find a point nearest (in some sense) to a cone and to spheres in the image of affine transformations.
		\vspace{-0.0in}}}
	\end{center}
There are two very prevalent examples of this problem that appear regularly in the literature, namely phase retrieval and source localization. There is no shortage of numerical schemes for addressing these problems, but most of these methods can be grouped into less than a handful of distinctly different ideas. We compare these different classes of algorithms, 13 algorithms in all (not counting equivalent algorithms under different names), on two different source localization problems (with and without noise) and six different phase retrieval problems, four synthetic and two involving laboratory data. A comparable benchmark study of (then) state of the art algorithms for phase retrieval was carried out by Marchesini \cite{Marchesini07} over a decade ago. In the intervening years, phase retrieval has been rediscovered by many in the statistics and applied mathematics communities. At the same time, there has been a steady development of the understanding of fixed point algorithms built on prox mappings for 
nonconvex problems. Our goal is to provide an update to the state of the art similar to \cite{Marchesini07} in the context of recent developments in nonconvex optimization.

The data and software for our experiments are all publicly accessible \cite{PTB}, and will hopefully serve as a starting point for future numerical methods. Amongst the algorithms we test, one has not appeared in the literature before, and this appears to be the best performer over all numerical experiments but one. The results expose some unwarranted enthusiasm for recently proposed approaches for phase retrieval, and show that algorithms that are beloved in the convex setting are not the best performers for nonconvex problems.

We are interested only in fixed points of algorithms and the relation of these points to the problem one {\em thinks} these algorithms solve. We present the algorithms through the lens of feasibility problems and their relaxations to smooth optimization. It is an unfortunate habit for researchers to conflate the fixed points of their algorithms with solutions to some optimization problem. Even when an algorithm is derived from a particular variational problem, the fixed points of the algorithm need not coincide with critical points of the motivating optimization problem, much less solutions. Our study of phase retrieval will hopefully remind readers of this distinction. Indeed, {\em no honest feasibility model of a phase retrieval problem has a solution}, yet, still, the algorithms converge to good fixed points that defy any obvious variational characterization. In fact, the algorithms built on feasibility models outperformed all other classes of algorithms by every relevant performance measure: 
computational efficiency, speed of convergence, accuracy and robustness against local minimums. It is often claimed that the analysis of algorithms for nonconvex feasibility problem is not well-understood, but this is not the case. A quantitative local analysis of algorithms for feasibility formulations of phase retrieval can be carried out with the tools and framework introduced in \cite{LukNguTam17}. A global analysis is incomplete, but some progress has been made in \cite{BST2018}. This is not to say that the book is closed for this problem, but at least the path forward is clear.

\section{Problem Instances}
We begin with a brief description of the two fundamental problem instances from which our numerical benchmarks are taken.

\subsection{Phase Retrieval}
The phase retrieval problem involves determining the phase of a complex valued function (that is, the real and imaginary parts) from the amplitude of its image under a unitary linear mapping - most often a Fourier transform, which is denoted by $\Fcal$. The problem is fundamental to diffraction imaging experiments. This includes near-field holography, far field diffraction imaging, ptychography, and wavefront sensing to name just a few observation techniques where this problem arises, see e.g., \cite{Luke17} and references therein. In all instances the measurements are real, nonnegative numbers which we denote by $b_i$ at the $i$-th pixel of some imaging array. A single image including all array elements is simply a vector $b \in \Rn$ where $n$ is the total number of pixels.

In many applications experimentalists often take several images at different settings of their instrument. It will be convenient in what follows to expand the measurement vector $b$ into an $n\times m$ real nonnegative matrix, each column of which represents a single image, $b_{ij}$ then being the $i$-th observation of the $j$-th image with $m$ images total. The different settings of the measurement device are accounted for in different operators $\Pcal_{j}$ -- the mapping $\Fcal$ at the foundation of the imaging model remains unchanged. There are a number of different experimental settings that can be captured by the mapping $\Pcal_{j}$; we mention a few here. In a typical ptychographic experiment an x-ray beam scans a specimen from side to side and an image is recorded at each scan position of the beam. The $j$-th spatial translation of the beam is represented by $\Pcal_{j}$. The difference here to the model discussed in \cite{HLST}, for instance, is that the beam is assumed to be fully known and shift 
invariant. Alternatively, as discussed in \cite{Luke02a, Hagemann14} one could record the image at different positions along the axis of propagation of the field, again represented by different mappings $\Pcal_{j}$ which account for magnification and defocus factors involved in such spatial translations. More recently, there have been proposals for adding random phase masks at the object plane. The idea of phase masks has its detractors \cite{Luke17}, but we include this data set as a point of reference to the more recent phase retrieval literature and it too fits the general model: the $m$ different masks can be represented by $\Pcal_{j}$ for $j = 1 , 2 ,  \ldots , m$.

In all the above cases, the approximate physical model for the measurement is
\begin{equation} \label{eq:physical model}
\left\|\paren{\Fcal \cdot \Pcal_{j}(z)}_{i}\right\| = b_{ij}, \quad \forall \,\, j = 1 , 2 , \ldots , m, \,\,\forall \,\, i = 1 , 2 , \ldots , n.
\end{equation}
In addition to these data equations, there are certain a priori qualitative constraints that can (and should) be added depending on the type of experiment that has been conducted. Often these are support constraints, or real-valuedness, or nonnegativity. In astronomy, for instance, the constraints are often support and magnitude constraints describing a field of constant magnitude across the aperture of the telescope. These types of constraints are described by the model \eqref{eq:physical model} where $\Fcal$ is the identity mapping. In x-ray diffraction imaging support or support and nonnegativity constraints are frequently imposed. In the early $00$'s Oszl\'anyi and S\"ut\'o \cite{Oszlanyi04,Oszlanyi08} proposed a simple {\em charge flipping} procedure that was quickly integrated into the software of the crystallography community for structure determination \cite{Smaalen03, PalCha07, Coelho07}.  Marchesini \cite{Marchesini08} identified the charge flipping operation as a reflector of a hard thresholding 
operation in the first algorithm for {\em sparse phase retrieval}. Sparsity constraints were more recently applied by Loock and Plonka \cite{LoockPlonka14,LoockPlonka16} for phase retrieval of objects enjoying a sparse representation in a wavelet-type dictionary.

The solution to the phase retrieval problem as presented here is the complex-valued vector $z \in \Cbb^{n}$ that satisfies \eqref{eq:physical model} for all $i$ and $j$ in addition to the a priori information implicit in the experiment (support, nonnegativity, sparsity, etc). Representing $n$-dimensional complex valued vectors instead as two-dimensional vectors on an $n$-dimensional product space, the phase retrieval problem is to find $z = \left(z_{1} ,  z_{2} , \ldots , z_{n}\right) \in \paren{\rr^{2}}^{n}$ ($z_{j}\in \rr^{2}$) satisfying \eqref{eq:physical model} for all $i$ and $j$ in addition to qualitative constraints. We return below to the issue of existence of such a vector. The short answer is that, unless the data is the pattern created by a periodic object, there does not exist a solution to problem \eqref{eq:physical model}, much less a {\em unique} solution. But this is not a concern for us, either practically or theoretically. One of the main goals of our study here is to convince readers that 
there are mathematically sound approaches that do not involve unrealistic assumptions of existence and uniqueness of solutions to equations that do not possess solutions.

The models and methods discussed here are based on a {\em feasibility} approach to this problem, that is at least initially, we will be happy with {\em any} point that comes as close as possible to matching the measurements. This, in short, is the way around existence and uniqueness. The points satisfying single measurements can be represented by sets. With this perspective, the question of existence is transformed into the question of whether the sets corresponding to different measurements have points in common. The question of uniqueness amounts to whether there is only a single common point between the sets.

\subsection{Source Localization}
The source localization problem appears in a broad range of applications including, for instance, mobile communication and wireless networks \cite{CSMC2004, SL2006}, acoustic/sound source localization \cite{LZBK2011}, GPS localization \cite{BP2012}, and brain activity identification \cite{F2010}, to mention just a few. The problem is based on distance measurements from an array of sensors (also called  anchors).  Here one is given a collection of $m$ sensors which are denoted by $a_{j} \in \rr^{d}$, where $d = 2$ or $3$ and $j = 1 , 2 , \ldots , m$. Each $a_{j}$ contains the exact location of the $j$-th sensor. The data consists of distance measurements between the unknown source and the sensors; this data is represented by $b_{j} > 0$, $j = 1 , 2 , \ldots , m$, and denotes the (possibly noisy) measurement of the range between the source and the $j$-th sensor $a_{j}$. This is described by the following equations:
\begin{equation} \label{Ranges}
	b_{j} = \norm{\Pcal_{j}(z)}, \quad j = 1 , 2 , \ldots , m,
\end{equation}
where $\Pcal_{j}(z) \equiv z - a_{j}$, $j = 1 , 2 , \ldots , m$, is the linear shift mapping. The problem is then to find an adequate approximation of the unknown source $\zbar$ satisfying the system \eqref{Ranges}. The only difference between source localization and the phase retrieval problems presented in the previous section is that source localization does not involve a Fourier-like transform $\Fcal$ (see model \eqref{eq:physical model}).

\subsection{Unifying Representation}
There are a number of different ways to represent the sets which we show are all equivalent with regard to the algorithms. We place all models in a real vector space, though in the context of phase retrieval it is understood that this is a reformulation of a complex valued model space. The model mappings in the previous sections are then linear mappings $\Fcal~:~ \paren{{\rr^{d}}}^{n} \to \paren{{\rr^{d}}}^{n}$ and $\Pcal_{j}~:~ \paren{{\rr^{d}}}^{n} \to \paren{{\rr^{d}}}^{n}$ with $\Fcal \cdot\Pcal_{j}(z) = \zhat = \left(\zhat_{1} , \zhat_{2} , \ldots , \zhat_{n}\right)$ for $\zhat_{i} \in\rr^{d}$. The sets of possible vectors satisfying the measurements are given by
\begin{equation}\label{eq:Cj}
	C_{j} \equiv \set{z \in \paren{\rr^{d}}^{n}}{\left\|\paren{\Fcal \cdot\Pcal_{j}(z)}_{i}\right\| = b_{ij}, \quad \forall \,\, i = 1 , 2 , \ldots , n}.
\end{equation}
Alternatively, one can work with the sets
\begin{equation}\label{eq:C'j}
	C_{j}' \equiv \set{z \in \paren{\rr^{d}}^{n}}{\left\|\paren{\Fcal(z)}_{i}\right\| = b_{ij}, \quad \forall \,\, i = 1 , 2 , \ldots , n},
\end{equation}
or
\begin{equation}\label{eq:Chatj}
	\Chat_{j} \equiv \set{z \in \paren{\rr^{d}}^{n}}{\left\|z_{i}\right\| = b_{ij}, \quad \forall \,\, i = 1 , 2 , \ldots , n}.
\end{equation}
Note that, for all $j = 1 , 2 , \ldots , m$, we have the following relations between these sets
\begin{equation}\label{eq:C_j relations}
	C_{j} = \Pcal_{j}^{\ast}C_{j}' = \Pcal_{j}^{\ast}\Fcal^{\ast}\Chat_{j},
\end{equation}
where $\Pcal_{j}^{\ast}$ and $\Fcal^{\ast}$ are the adjoints of $\Pcal_{j}$ and $\Fcal$, respectively.

For phase retrieval $d = 2$, $n$ is large ($1024^{2}$ is not uncommon) and $m$ is anywhere from $1$ to $10$. For sensor localization, $d = 2$ or $3$, $n = 1$ and $m$ is at least $3$, but not usually greater than $100$. These sets are in different spaces relative to one another, but for both applications $\Fcal$ and $\Pcal_{j}$, $j = 1 , 2 , \ldots , m$, are often unitary, so the choice of space to work in is a matter of convenience only. Even though they are nonconvex (a line segment between any two points in the sets does not belong to the sets), phase sets have very nice structure. The sets $\Chat_{j}$, for instance, are just $\ell_{2}$-spheres in an $n$-dimensional product space of $\paren{\rr^{d}}$ with component-wise radii given by the elements of the vector $b$. This is true regardless of whether or not the measurement $b$ is contaminated with noise. As such, these sets are smooth, {\em semialgebraic} (constructible by finite systems of polynomial inequalities), and {\em prox-regular} (loosely defined 
as sets with locally single-valued projections \cite{PolRockThib00}). By the relationship \eqref{eq:C_j relations}, the sets $C_{j}$ and $C_{j}'$ also enjoy this nice structure. These facts were already observed in \cite[Proposition 3.5]{HLST} and \cite{Luke12}.

We reserve the set $C_{0}$ for the qualitative constraints. The qualitative constraints that most often applied are either of the same form as the sets above - as in the case of support and magnitude constraints - or cones. A support constraint alone is a restriction of points to a subspace. Support and nonnegativity is the positive orthant of a subspace, that is, a convex cone. It is easy to see that sets of points satisfying sparsity constraints can be characterized as unions of subspaces (see, for instance, \cite{BauLukePhanWang14}), hence the set of points satisfying a sparsity constraint is also a cone.

\section{Variational Models}
A key feature that allows one to classify algorithms is the degree of smoothness in the underlying model.  We group the algorithms into two classes below in increasing order of model smoothness. A third class of algorithms mixes the first two, and involves a {\em product space formulation} that gathers constraints/sets into blocks at a modest cost in increased dimensionality. Nonsmooth aspects are maintained within the blocks, but communication between blocks is modeled, most often, by smooth operators. The numerical results presented in Section \ref{s:numerics} indicate that the smoother the model is, the slower the algorithm progresses toward a fixed point. Smoothness is also sometimes associated with stability or reliability. Here again, our results do not support this intuition for this family of nonconvex problems:  the most reliable and robust algorithms are also in the nonsmooth and feasibility-based category. We do not take into account advantages of parallelization and other architecture-dependent 
features that could make a difference in clock times.

\subsection{Model Category I: Multi-Set Feasibility} \label{SS:Feas}
The most natural place to start is by naively trying to find a point in the intersection of the data generated sets $C_{j}$, $j = 1 , 2 , \ldots , m$, given by \eqref{eq:Cj} together with the possible qualitative constraint set $C_{0}$:
\begin{equation*}
	\Find z^{\ast} \in \bigcap_{j = 0}^{m} C_{j}.
\end{equation*}
We show in Section \ref{s:numerics} that this leads to the most effective methods for solving the problem.  Keeping in mind, however, that for all practical purposes the intersection above is empty, we also show that the algorithm is not solving the problem we thought it was solving.

For our purposes we will prefer to pose this problem equivalently in an optimization format
\begin{equation}\label{eq:sum-of-indicators}
	\min_{z \in \paren{\rr^{d}}^{n}}\sum_{j = 0}^{m} \iota_{C_{j}}\left(z\right),
\end{equation}
where
\begin{equation*}
	\iota_{C_{j}}\left(z\right) \equiv
	\begin{cases}
		0, & \mbox{ if } z \in C_{j}, \\
    		+\infty, & \mbox{ else}.
	\end{cases}
\end{equation*}
The indicator function $\iota_{C_{j}}$, $j = 0 , 1 , \ldots , m$, is an {\em extended real-valued} function from $\paren{\rr^{d}}^{n}$ to the extended real-line $\rr \cup \{ +\infty \}$. The fact that the intersection is empty is reflected in the fact that the optimal value to problem \eqref{eq:sum-of-indicators} is $+\infty$.

Despite these worrisome issues, we examine the {\em Method of Cyclic Projections}:

    \begin{center}
    		\fcolorbox{black}{Ivory2}{\parbox{15.7cm}{\vspace{-0.1in}
		\begin{alg}[{\bf Cyclic Projections - CP}]\label{alg:cp} $~$\\
    			{\bf Initialization.} Choose $z^{0} \in \paren{\rr^{d}}^{n}$. \\
    			{\bf General Step ($k = 0 , 1 , \ldots$)}
			    \begin{equation*}
					z^{k + 1} \in P_{C_{0}}P_{C_{1}} \cdots P_{C_{m}}z^{k}.
			    \end{equation*}
		\end{alg}\vspace{-0.1in}}}
	\end{center}

Here the {\em orthogonal projector} $P_{C_{j}}$ for $j = 0 , 1 , \ldots , m$ is defined by
\begin{equation*}
P_{C_{j}}\left(z\right) \equiv \argmin_{y \in C_{j}} \norm{y - z}.
\end{equation*}
Since $C_{j}$ is nonconvex the projector is, in general, a set-valued mapping. This can be computed explicitly \cite[Corollary 4.3]{Luke02a}, for all $j = 1 , 2 , \ldots , m$, by the formula
\begin{equation}\label{eq:P_Cj}
	y \in P_{C_{j}}\left(z\right) \quad \iff y = \Pcal_{j}^{\ast}\Fcal^{\ast}\yhat, \quad \yhat_{i} \in
	\begin{cases}
		b_{ij}\frac{\paren{\Fcal\Pcal_{j}(z)}_{i}}{\left\|\paren{\Fcal\Pcal_{j}(z)}_{i}\right\|}, & \mbox{ if }\paren{\Fcal\Pcal_{j}(z)}_{i} \neq 0, \\
		b_{ij}\Sbb, & \mbox{ if } \paren{\Fcal\Pcal_{j}(z)}_{i} = 0.
	\end{cases}
\end{equation}
The unit sphere in $\rr^{d}$ is denoted by $\Sbb$ above. The projector $P_{C_{0}}$ is also explicitly known and has a structure that is no more complicated than \eqref{eq:P_Cj}, often it is simpler.

The analysis of this algorithm for consistent and inconsistent nonconvex problems has been established in \cite{LukNguTam17}. In the inconsistent case, the algorithm fixed points need not correspond to best approximation points. Instead, the fixed points generate cycles of smallest length locally over all other possible cycles generated by projecting onto the sets in the same order. Strong guarantees on convergence of the iterates, like a local linear rate, have been established generically for problems with this structure (see \cite[Example 3.6]{LukNguTam17}). What remains is to establish guarantees for global convergence to fixed points.

Another cyclic projection-type method for inconsistent feasibility problems is based on what is commonly known as the {\em Douglas-Rachford (DR)} Algorithm (the original algorithm is a domain splitting method for partial differential equations \cite{DougRach56}). In our context the method can only be applied directly to two-set feasibility problems,
\begin{equation*}
	\Find x \in C_{0} \cap C_{1}.
\end{equation*}
The fixed point iteration is given by
\begin{equation}\label{eq:dr}
	(DR) \qquad\qquad\qquad z^{k + 1} \in \frac12 \paren{R_{C_{0}}R_{C_{1}} + \Id}z^{k},
\end{equation}
where $R_{C} \equiv 2P_{C} - \Id$ is the \textit{reflector} operator of the set $C$ and $\Id$ denotes the identity operator. It is important not to forget that, even if the feasibility problem is consistent, the fixed points of the Douglas-Rachford Algorithm will not in general be points of intersection. Instead, the {\em shadows} of the iterates defined as $P_{C_{1}}z^{k}$, $k \in \nn$, converge to intersection points, when these exist \cite{BCL3}.

To extend this to more than two sets, Borwein and Tam \cite{borwein2014cyclic, BorweinTam15} proposed the following variant:

    \begin{center}
    		\fcolorbox{black}{Ivory2}{\parbox{15.7cm}{\vspace{-0.1in}
		\begin{alg}[{\bf Cyclic Douglas-Rachford - CDR}]\label{alg:cdr} $~$\\
    			{\bf Initialization.} Choose $z^{0} \in \paren{\rr^{d}}^{n}$. \\
    			{\bf General Step ($k = 0 , 1 , \ldots$)}
			      \begin{equation*}
				    z^{k + 1} \in \paren{\frac12 \paren{R_{C_{0}}R_{C_{1}} + \Id}}\paren{\frac12 \paren{R_{C_{1}}R_{C_{2}} + \Id}}\cdots\paren{\frac12 \paren{R_{C_{m}}R_{C_{0}} + \Id}}z^{k}.
			      \end{equation*}
		\end{alg}\vspace{-0.1in}}}
	\end{center}

It is easy to envision different sequencing strategies than the one presented above. In \cite{BauschkeNollPhan15} one of the pair of sets is held fixed, and this has some theoretical advantages in a convex setting. We did not observe any advantage for the problems studied here. A comprehensive investigation of optimal sequencing strategies for problems with different structures has not been carried out.

A relaxation to the Douglas-Rachford Algorithm first proposed in \cite{Luke05a} is described below in Algorithm \ref{alg:raar}. At this stage, we just motivate it as a convex combination of the Douglas-Rachford mapping and the projection onto the ``inner'' set $C_{1}$ above: for $\lambda \in \left(0 , 1\right]$
\begin{equation} \label{eq:raar}
	(DR\lambda) \qquad\qquad\qquad z^{k + 1} \in \paren{\frac{\lambda}{2}  \paren{R_{C_{0}}R_{C_{1}} + \Id} + \left(1 - \lambda\right)P_{C_{1}}}z^{k}.
\end{equation}
Extending this to more than two sets yields the following algorithm, which has not appeared in the literature before.

    \begin{center}
    		\fcolorbox{black}{Ivory2}{\parbox{15.7cm}{\vspace{-0.1in}
		\begin{alg}[{\bf Cyclic Relaxed Douglas-Rachford - CDR$\lambda$}]\label{alg:caarl} $~$\\
    			{\bf Initialization.} Choose $z^{0} \in \paren{\rr^{d}}^{n}$ and $\lambda \in \left[0 , 1\right]$. \\
    			{\bf General Step ($k = 0 , 1 , \ldots$)}
				\begin{eqnarray*}
					z^{k + 1} \in & \paren{\frac{\lambda}{2} \paren{R_{C_{0}}R_{C_{1}} + \Id} + \left(1 -\lambda\right)P_{C_{1}}}\paren{\frac{\lambda}{2} \paren{R_{C_{1}}R_{C_{2}} + \Id} + \left(1 - \lambda\right)P_{C_{2}}} & \\
					& \cdots \paren{\frac{\lambda}{2} \paren{R_{C_{m}}R_{C_{0}} + \Id} + \left(1 -\lambda\right)P_{C_{0}}}z^{k}.
			      \end{eqnarray*}
		\end{alg}\vspace{-0.1in}}}
	\end{center}

The analysis for DR$\lambda$ and its precursor, Douglas-Rachford is contained in \cite[Section 3.2.2]{LukNguTam17}. The crucial property to be determined is any type of metric subregularity of the fixed point mapping at its fixed points. This remains open for practical phase retrieval, though we believe that generic guarantees of {\em some} sort of convergence (most generically, sublinear) are readily obtainable from current results.

In the convex setting the Douglas-Rachford Algorithm can be derived from the Alternating Directions Method of Multipliers (ADMM, \cite{Glowinski75}) for solving a dual problem \cite{Gabay83}. This algorithm is extremely popular at the moment for large-scale convex problems with linear constraints, and the literature in this context is massive (see, for instance, \cite{ST2013} and references therein). In the nonconvex setting, the dual correspondence is lost, though there have been some recent developments and studies \cite{LP2015,Patrinos17,BST2018}. The ADMM falls into the category of augmented Lagrangian-based methods. Thus, we can reformulate problem \eqref{eq:sum-of-indicators} as 
\begin{equation}\label{eq:sum-of-indicators2}
	\min_{x , z_{j} \in \paren{\rr^{d}}^{n}} \set{\iota_{C_{0}}\left(z\right) + \sum_{j = 1}^{m} \iota_{C_{j}}\left(z_{j}\right)}{z_{j} = x, \,\, j = 1 , 2 , \ldots , m},
\end{equation}
and then one can apply ADMM to the augmented Lagrangian of this formulation; specifically, we have
\begin{equation}\label{eq:aug Lagr1}
	\Ltilde_{\eta}\left(x , z_{j} , v_{j}\right) \equiv \iota_{C_{0}}\left(x\right) + \sum_{j = 1}^{m} \left(\iota_{C_{j}}\left(z_{j}\right) + \ip{v_{j}}{x - z_{j}} + \frac{\eta}{2}\norm{x - z_{j}}^{2}\right),
\end{equation}
where $\eta > 0$ is a penalization parameter and $v_{j}$, $j = 1 , 2 , \ldots , m$, are the multipliers which associated with the linear constraints. The ADMM algorithm applied to finding the critical points of the corresponding augmented Lagrangian (see \eqref{eq:aug Lagr1}) is given by

    \begin{center}
    		\fcolorbox{black}{Ivory2}{\parbox{15.7cm}{\vspace{-0.1in}
		\begin{alg}[{\bf Nonsmooth ADMM$_1$}]\label{alg:ADMM1} $~$\\
    			{\bf Initialization.} Choose $x^{0} , z_{j}^{0} , v_{j}^{0} \in \paren{\rr^{d}}^{n}$ and fix $\eta>0$. \\
    			{\bf General Step ($k = 0 , 1 , \ldots$)}
			\begin{itemize}
				\item[1.] Update
					\begin{align} \label{ADMM:Step1}
						x^{k + 1} & \in \argmin_{x \in \paren{\rr^{d}}^{n}} \left\{ \iota_{C_{0}}\left(x\right) + \sum_{j = 1}^{m} \left(\ip{v_{j}^{k}}{x - z_{j}^{k}} + \frac{\eta}{2}\norm{x - z_{j}^{k}}^{2}\right) \right\} \nonumber \\
						& = P_{C_{0}}\left(\frac{1}{m}\sum_{j = 1}^{m} \left(z_{j}^{k} - \frac{1}{\eta}v_{j}^{k}\right)\right).
					\end{align}
				\item[2.] For all $j = 1 , 2 , \ldots , m$ update (in parallel)
					\begin{align} \label{ADMM:Step2}
						z_{j}^{k + 1} & \in \argmin_{z_{j} \in \paren{\rr^{d}}^{n}} \left\{ \iota_{C_{j}}\left(z_{j}\right) + \ip{v_{j}^{k}}{x^{k + 1} - z_{j}} + \frac{\eta}{2}\norm{x^{k + 1} - z_{j}}^{2} \right\} \nonumber \\
					& = P_{C_{j}}\left(x^{k + 1} - \eta v_{j}^{k}\right).
					\end{align}
				\item[3.] For all $j = 1 , 2 , \ldots , m$ update (in parallel)
					\begin{equation} \label{ADMM:Step3}
						v_{j}^{k + 1} = v_{j}^{k} + \eta\paren{x^{k + 1} - z_{j}^{k + 1}}.		
					\end{equation}
			\end{itemize}
		\end{alg}\vspace{-0.2in}}}
	\end{center}

An ADMM scheme for phase retrieval has appeared in \cite{LSJL2018} and for sensor localization in \cite{LukSabTebZat17}. Our numerical experiments below do not indicate any advantage of this approach. We include it, however, as a point of reference to both the Douglas-Rachford Algorithm and the smoother ADMM$_2$ Algorithm discussed below in Algorithm \ref{alg:admm-locss}. Seeing what changes render a bad algorithm reasonable  provides tremendous insight. Note that the projections in Step 2 of the algorithm can be computed in parallel, while the Cyclic Projections and Cyclic Douglas-Rachford Algorithms must be executed sequentially. The benchmarking comparisons carried out in Section \ref{s:numerics} do not reflect the advantages of parallelizable methods like ADMM$_1$ when implemented on multiple CPU's/GPU's architectures.

\subsection{Model Category II: Smooth Nonconvex Optimization}
The next algorithm, {\em Averaged Projections}, could be motivated purely from the feasibility framework detailed above. However, since there is a more significant {\em smooth} interpretation of this model, we present it here as a gateway to the smooth model class.

    \begin{center}
    		\fcolorbox{black}{Ivory2}{\parbox{15.7cm}{\vspace{-0.1in}
		\begin{alg}[{\bf Averaged Projections - AvP}]\label{alg:avp} $~$\\
    			{\bf Initialization.} Choose $z^{0} \in \paren{\rr^{d}}^{n}$. \\
    			{\bf General Step ($k = 0 , 1 , \ldots$)}
			    \begin{equation*}
					z^{k + 1} \in \frac{1}{m + 1}\sum_{j = 0}^{m} P_{C_{j}}z^{k}.
			    \end{equation*}
		\end{alg}\vspace{-0.1in}}}
	\end{center}

This algorithm is often preferred not only because it is parallelizable, but also because it appears to be more robust to problem inconsistency. Indeed, Averaged Projections Algorithm can be equivalently viewed as gradient-based schemes when applied to an adequate smooth and nonconvex objective function. This well-known fact goes back to \cite{Zarantonello} when the sets $C_{j}$, $j = 0 , 1 , \ldots , m$, are closed and convex. To be specific, we describe below two such schemes which are shown to be equivalent to AvP (see also Section 3.3 for further equivalent views), followed by a Dynamically Reweighted Averaged Projections Algorithm.

First, consider the problem of minimizing the sum of squared distances to the sets $C_{j}$, $j = 0 , 1 , \ldots , m$, that is,
\begin{equation}\label{eq:sum-of-sq-dist}
	\ucmin{f\left(z\right) \equiv \frac{1}{2\left(m + 1\right)}\sum_{j = 0}^{m} \dist^{2}\left(z , C_{j}\right)}{z \in \paren{\rr^{d}}^{n}}.
\end{equation}
Since the sets $C_{j}$, $j = 0 , 1 , \ldots , m$, are nonconvex, the functions $\dist^{2}\left(z , C_{j}\right)$ are clearly not differentiable, and hence, same for the objective function $f\left(z\right)$. However, in our context, the sets $C_{j}$, $j = 0 , 1 , \ldots , m$, are prox-regular (cf. Section 2.3). From elementary properties of prox-regular sets \cite{PolRockThib00} it can be shown that the gradient of the squared distance is defined and differentiable  with Lipschitz continuous derivative (that is, the corresponding Hessian)  up to the boundary of $C_{j}$, $j = 0 , 1 , \ldots , m$, and points where the coordinate elements of the vector $z$ vanish. Indeed, for $f$ given by \eqref{eq:sum-of-sq-dist} we have
\begin{equation}\label{eq:grad sum-of-sq-dist}
	\nabla f\left(z\right) \equiv \frac{1}{m+1}\sum_{j = 0}^{m} \paren{\Id - P_{C_{j}}}\left(z\right).
\end{equation}
Thus, applying the gradient descent with {\em unit stepsize} to problem \eqref{eq:sum-of-sq-dist}, one immediately recovers the AvP.

Readers familiar with variational analysis will also recognize \eqref{eq:sum-of-sq-dist} as the relaxation of \eqref{eq:sum-of-indicators} via the {\em Moreau envelope} \cite{Moreau65} of the indicator functions $\iota_{C_{j}}$, $j = 0 , 1 , \ldots , m$. But even without these sophisticated tools, inspection shows that the objective in \eqref{eq:sum-of-sq-dist} is smooth, has full domain and takes the value zero at points of intersection, if such points exist. These kinds of models are much more prevalent in applications than the more severe-looking feasibility format of problem \eqref{eq:sum-of-indicators}.

There is another interesting way to tackle problem \eqref{eq:sum-of-sq-dist}, allowing one to make links with another fundamental algorithmic paradigm. Ignoring the weighting factor for the moment, we consider the following problem:
\begin{equation} \label{SumModel:1}
	\min_{z \in \paren{\rr^{d}}^{n}} f\left(z\right) \equiv \frac{1}{2}\sum_{j = 0}^{m} \dist^{2}\left(z , C_{j}\right).
\end{equation}
Using the definition of the function $\dist^{2}\left(\cdot , C_{j}\right)$, $j = 0 , 1 , \ldots , m$, we can reformulate problem \eqref{SumModel:1} as follows
\begin{equation} \label{SumModel:2}
	\min_{z , \bu} \set{\sum_{j = 0}^{m} \frac{1}{2}\norm{z - u_{j}}^{2}}{u_{j} \in C_{j}, \quad j = 0 , 1 , \ldots , m},
\end{equation}
where $\bu = \left(u_{0} , u_{1} , \ldots , u_{m}\right) \in \left(\paren{\rr^{d}}^{n}\right)^{m + 1}$.

Note that problem \eqref{SumModel:2} always has an optimal solution (since we minimize continuous function over a closed and bounded set).
The structure of the optimization problem \eqref{SumModel:2}, which includes constraint sets that is separable over the variables $u_{j}$, $j = 0 , 1 , \ldots , m$, suggests that one can exploit this generous property when developing an optimization algorithm. Alternating Minimization (AM) is a classical optimization technique which was designed exactly for these situations, and involves updating each variable separately in a cyclic manner. More precisely, AM when applied to problem \eqref{SumModel:2}, generates sequences defined by the following algorithm.

    \begin{center}
    		\fcolorbox{black}{Ivory2}{\parbox{15.7cm}{\vspace{-0.1in}
		\begin{alg}[{\bf Alternating Minimization - AM}]\label{alg:am} $~$\\
    			{\bf Initialization.} Choose $\left(z^{0} , u_{0}^{0} , u_{1}^{0} , \ldots , u_{m}^{0}\right) \in \left(\paren{\rr^{d}}^{n}\right)^{m + 2}$. \\
    			{\bf General Step ($k = 0 , 1 , \ldots$)}
			\begin{itemize}
				\item[1.] Update
					\begin{equation} \label{AM:Step1}
						z^{k + 1} = \argmin_{z \in \paren{\rr^{d}}^{n}} \sum_{j = 0}^{m} \frac{1}{2}\norm{z - u_{j}^{k}}^{2} =  \frac{1}{m + 1} \sum_{j = 0}^{m} u_{j}^{k}.
					\end{equation}
				\item[2.] For all $j = 0 , 1 , \ldots , m$ update (in parallel)
					\begin{equation} \label{AM:Step2}
						u_{j}^{k + 1} \in \argmin_{u_{j} \in C_{j}} \frac{1}{2}\norm{u_{j} - z^{k + 1}}^{2} = P_{C_{j}}\left(z^{k + 1}\right).
					\end{equation}
			\end{itemize}
		\end{alg}\vspace{-0.2in}}}
	\end{center}

This means that, by combining \eqref{AM:Step1} and \eqref{AM:Step2}, the algorithm can be compactly written as
\begin{equation} \label{AM:Equiv}
	z^{k + 1} \in \frac{1}{m + 1} \sum_{j = 0}^{m} P_{C_{j}}z^{k},
\end{equation}
which is exactly the Averaged Projections Algorithm \ref{alg:avp}. In the case that $m = 1$, i.e., only one image is considered, the Alternating Minimization Algorithm discussed above coincides with what was called the {\em Error Reduction} Algorithm in \cite{Fien78}. In \cite[Remark 2.2(i)]{HLST} the more general PHeBIE Algorithm applied to the problem of blind ptychography was shown to reduce to Averaged Projections Algorithm for phase retrieval when the illuminating field is known. The PHeBIE Algorithm is a slight extension of the PALM Algorithm \cite{BST2014}. A partially preconditioned version of PALM was studied in
\cite{Marchesini18} for phase retrieval, with improved performance over PALM.

The analysis of Averaged Projections Algorithm for problems with this geometry is covered by the analysis of nonlinear/nonconvex gradient descent. Much of this is classical and can be found throughout the literature, but it is limited to guarantees of convergence to {\em critical points} (see, for instance, \cite{ABS2013,BST2014}). This begs the question as to which of the critical points are global minima. The answer to this is unknown.

Instead of always taking the fixed average $1/(m + 1)$, in the formulation of problem \eqref{eq:sum-of-sq-dist}, it is possible to derive a variational interpretation of dynamically weighted averages between the projections to the sets $C_{j}$, $j = 0 , 1 , \ldots , m$. This idea was proposed in \cite{Luke02a} where it is called {\em extended least squares}. A similar approach was also proposed in \cite{BTC2008} where the resulting algorithm is called the Sequential Weighted Least Squares (SWLS) Algorithm. The underlying model in \cite{Luke02a} is the negative log-likelihood measure of the sum of squared set distances:
\begin{equation}\label{e:llssd}
	\ucmin{\sum_{j = 0}^{m} \ln\left(\dist^{2}\left(z , C_{j}\right) + c\right)}{z \in \paren{\rr^{d}}^{n}}, \qquad (c > 0).
\end{equation}
Gradient descent applied to this objective yields the following {\em Dynamically Reweighted Averaged Projections} Algorithm.

    \begin{center}
    		\fcolorbox{black}{Ivory2}{\parbox{15.7cm}{\vspace{-0.1in}
		\begin{alg}[{\bf Dynamically Reweighted Averaged Projections - DyRePr}]\label{alg:DyRePr} $~$\\
    			{\bf Initialization.} Choose $z^{0} \in \paren{\rr^{d}}^{n}$ and $c > 0$. \\
    			{\bf General Step ($k = 0 , 1 , \ldots$)}
					\begin{equation} \label{DyRePr:Step1}
						z^{k + 1} \in z^{k} -  \sum_{j = 0}^{m} \frac{2}{\left(\dist^{2}\left(z^{k} , C_{j}\right) + c\right)}\left(z^{k} - P_{C_{j}}\left(z^{k}\right)\right).
					\end{equation}
		\end{alg}\vspace{-0.2in}}}
	\end{center}

{\bf Beyond first order methods.} The twice continuous differentiability of the sum of squared distances at most points suggests that one could try higher-order techniques from nonlinear optimization in order to accelerate the basic gradient descent method. Higher-order accelerations, like quasi-Newton methods, rely on extra smoothness in the objective function. However, as observed in \cite{LewisOverton13}, quasi-Newton based method can be used to solve nonsmooth problems. In the numerical comparisons in Section \ref{s:numerics} we benchmark a limited memory BFGS method applied in \cite{Luke02a} against the other techniques.

    \begin{center}
    		\fcolorbox{black}{Ivory2}{\parbox{15.7cm}{\vspace{-0.1in}
		\begin{alg}[{\bf Limited Memory BFGS with Trust Region - QNAvP}]\label{alg:LBFGS} $~$\\ \vspace{-0.3in}
		\begin{enumerate}[1.]
			\item (Initialization) Choose $\tilde{\eta} > 0,$ $\zeta > 0,$ $\overline{\ell} \in \{ 1 , 2 ,\ldots , n \}$, $z^{0} \in \Cbb^{n}$, and set $\nu = \ell = 0$. Compute $\nabla f\left(z^{0}\right)$ and $\norm{\nabla f\left(z^{0}\right)}$ for
				\begin{equation*}
					f\left(z\right) \equiv \frac{1}{2\left(m + 1\right)}\sum_{j = 0}^{m} \dist^{2}\left(z , C_{j}\right), \quad \nabla f\left(z\right) \equiv \frac{1}{m + 1}\sum_{j = 0}^{m} \paren{\Id - P_{C_{j}}}\left(z\right).
				\end{equation*}
			\item (L-BFGS step) For each $k = 0 , 1 , 2 , \ldots$ if $\ell = 0$ compute $z^{k + 1}$ by some line search algorithm; otherwise compute
				\begin{equation*}
					s^{k} = -\left(M^{k}\right)^{-1}\nabla f\left(z^{k}\right),
		    		\end{equation*}
		    		where $M^{k}$ is the L-BFGS update \cite{ByrdNocedalSchnabel94}, $z^{k + 1} = z^{k} + s^{k}$, $f\left(z^{k + 1}\right)$, and the {\em predicted change} (see, for instance \cite{NocedalWright}).
		  	\item (Trust Region) If $\rho\left(s^{k}\right) < \tilde{\eta}$, where
		    		\begin{equation*}
		    			\rho\left(s^{k}\right) = \frac{\mbox{ actual change at step $k$}}{\mbox{ predicted change at step $k$}},
		    		\end{equation*}
		      	reduce the trust region $\Delta^{k}$, solve the trust region subproblem for a new step $s^{k}$ \cite{Burke1}, and return to the beginning of Step 2. If $\rho\left(s^{k}\right) \geq \tilde{\eta}$ compute $z^{k + 1} = z^{k} + s^{k}$ and $f\left(z^{k + 1}\right)$.
			\item (Update) Compute $\nabla f\left(z^{k + 1}\right)$, $\norm{\nabla f\left(z^{k + 1}\right)}$,
		      	\begin{equation*}
			  		y^{k} \equiv \nabla f\left(z^{k + 1}\right) - \nabla f\left(z^{k}\right), \quad s^{k}\equiv z^{k + 1} - z^{k},
		      	\end{equation*}
		    		and ${s^{k}}^{T}y^{k}$. If ${s^{k}}^{T}y^{k} \leq \zeta$, discard the vector pair $\{s^{k -\ell} , y^{k - \ell}\}$ from storage, set $\ell = \max\{\ ell - 1 , 0 \}$, $\Delta^{k + 1} = \infty$, $\mu^{k + 1} = \mu^{k}$ and $M^{k + 1} = M^{k}$ (\ie shrink the memory and don't update); otherwise set $\mu^{k + 1} = \frac{{y^{k}}^{T}y^{k}}{{s^{k}}^{T}y^{k}}$ and $\Delta^{k + 1} = \infty,$ add the vector pair $\{ s^{k} , y^{k} \}$ to storage, if $\ell = \overline{\ell}$, discard the vector pair $\{ s^{k - \ell} , y^{k -\ell} \}$ from storage. Update the Hessian approximation $M^{k + 1}$ \cite{ByrdNocedalSchnabel94}. Set $\ell = \min\{ \ell + 1 , \overline{\ell} \}$, $\nu = \nu + 1$ and return to Step 1.
    		\end{enumerate}
		\end{alg}\vspace{-0.1in}}}
	\end{center}
	
This rather complicated-looking algorithm is a standard in nonlinear optimization and it even shows unexpected (and largely unexplained) good performance for nonsmooth problems \cite{LewisOverton13}. Convergence is still open, but we include this algorithm as one of the typical types of accelerations via higher order methods one might try.

{\bf Least squares based-models.} In \cite{Marchesini07b} Marchesini studies an augmented Lagrangian approach to solving
\begin{equation} \label{LS}
	\min_{z \in \paren{\rr^{d}}^{n}}\frac{1}{2n}\sum_{j = 0}^{m}\sum_{i = 1}^{n} \left(\left\|\paren{\Fcal \cdot \Pcal_{j}(z)}_{i}\right\| - b_{ij}\right)^{2}.
\end{equation}
It is not difficult to see that this is a nonsmooth least-squares relaxation of problem \eqref{eq:physical model} \cite[Lemma 3.1]{PauBecEldSab18}. In reduced form, the resulting primal-dual/ADMM Algorithm has the following prescription which corrects an error in \cite{LukSabTebZat17}\footnote{The simplified form of \cite[Algorithm 1 (simplified)]{LukSabTebZat17} has an error in equation (5.2).}.

	\begin{center}
 		\fcolorbox{black}{Ivory2}{\parbox{15.7cm}{\vspace{-0.1in}
		\begin{alg}[ADMM$_2$] \label{alg:admm-locss} $~$\\
    			{\bf Initialization.}  Choose any $x^{0} \in \rr^{n}$ and $\rho_{j} > 0$, $j = 0 ,  1 , \ldots , m$. Compute $u_{j}^{1} \in P_{C_{j}}\left(z^{0}\right)$ ($j = 0 , 1 , \ldots , m)$ and $z^{1} \equiv \left(1/\left(m + 1\right)\right)\sum_{j = 0}^{m} u_{j}^{1}$. \\
			{\bf General Step.} For each $k = 1 , 2 , \ldots$ generate the sequence $\left\{ \left(z^{k} ,
				\bu^{k}\right) \right\}_{k \in \nn}$ as follows:
            		\begin{itemize}
                		\item Compute
                    		\begin{equation} \label{ADM:StepXs}
                        		z^{k + 1} = \frac{1}{m}\sum_{j = 1}^{m} \left(u_{j}^{k} + \frac{1}{\rho_{j}}\left(z^{k} - z^{k - 1}\right)\right).
                    		\end{equation}
                		\item For each $j = 1 , 2 , \ldots , m$, compute
                   		\begin{equation} \label{ADM:StepUs}
				  			u_{j}^{k + 1} = P_{C_{j}}\left(u_{j}^{k} + \frac{1}{\rho_{j}}\left(2z^{k} - z^{k - 1}\right)\right).
                      	\end{equation}
				\end{itemize}
		  \end{alg} \vspace{-0.2in}}}
	\end{center}
\noindent This algorithm also can be viewed as a smoothed/relaxed version of Algorithm \ref{alg:ADMM1}.

The appearance of the norm in \eqref{LS} makes the analysis of the least squares approach inconvenient without the tools of nonsmooth analysis. A popular way around this, is to formulate \eqref{eq:physical model} as a system of quadratic equations:
\begin{equation} \label{eq:squared physical model}
	\left\|\paren{\Fcal \cdot \Pcal_{j}(z)}_{i}\right\|^{2} = b_{ij}^{2}, \quad \forall \,\, j = 1 , 2 , \ldots , m, \,\,\forall \,\, i = 1 , 2 , \ldots , n.
\end{equation}
The corresponding squared least squares residual of the quadratic model which is exceedingly smooth:
\begin{equation} \label{SLS}
	\min_{z \in \paren{\rr^{d}}^{n}} G\left(z\right) \equiv \frac{1}{2}\sum_{j = 0}^{m} \sum_{i = 1}^{n} \left(\left\|\paren{\Fcal \cdot \Pcal_{j}(z)}_{i}\right\|^{2} - b_{ij}^{2}\right)^{2}.
\end{equation}

A popular approach in the applied mathematics and statistics communities for solving the squared least squares formulation \eqref{SLS} is based on a trick from conic programming for turning quadratics into linear functions in a lifted space of much higher dimension. The lifted linear objective is still constrained in rank, which is a nonconvex constraint, but this is typically relaxed to a convex constraint. The idea, called {\em phase lift} when applied to phase retrieval \cite{CandesEldarStrohmerVononinski}, is not an efficient approach due to a number of reasons, not the least of which being the increase in dimension (the square of the dimensionality of the original problem).  Indeed, the phase lift method is not implementable on standard consumer-grade architectures for almost all of the benchmarking experiments conducted below.

Undeterred by the experience with phase lift, the recent paper \cite{candes2014phase} has inspired many studies of the Wirtinger Flow (WF) Algorithm for solving problem \eqref{SLS}. The WF method is a gradient descent algorithm applied in the image space of the mapping $\Fcal$ to minimize the function $G$.	This leads to the following method.

    \begin{center}
    		\fcolorbox{black}{Ivory2}{\parbox{15.7cm}{\vspace{-0.1in}
		\begin{alg}[{\bf Wirtinger Flow - WF}]\label{alg:Wirt} $~$\\
    			{\bf Initialization.} Choose $z^{0} \in \paren{\rr^{d}}^{n}$ and step-size $\mu > 0$. \\
    			{\bf General Step ($k = 0 , 1 , \ldots$)}
				\begin{equation} \label{FW}
					z^{k + 1} = z^{k} - \frac{\mu}{\norm{z^{0}}^{2}}\left(\left\|\paren{\Fcal \cdot \Pcal_{j}(z^{k})}_{i}\right\|^{2} - b_{ij}^{2}\right)z^{k}.
				\end{equation}
		\end{alg}\vspace{-0.2in}}}
	\end{center}

While smoothness makes the analysis nicer, the quartic objective has almost no curvature around critical points, which makes convergence of first order methods much slower than first order methods applied to nonsmooth objectives. See \cite[Section 5.2]{Luke02a} for a discussion of this. The numerical comparisons here support this conclusion.
	
\subsection{Model Category III: Product Space Formulations}
This third category of algorithms synthesizes, in some sense, the previous two. Here, the basic idea is to lift the problem to the product space $\left(\paren{\rr^{d}}^{n}\right)^{m + 1}$ which can be then formulated as a two-set feasibility problem
\begin{equation*}
	\Find \bz^{\ast} \in C \cap D,
\end{equation*}
where $\bz^{\ast} = \left(z_{0}^{\ast} , z_{1}^{\ast} , \ldots , z_{m}^{\ast}\right)$, $C := C_{0} \times C_{1} \times \cdots \times C_{m}$ and $D$ is the diagonal set of $\paren{\Rbb^{d}}^{n(m + 1)}$ which is defined by $\left\{\bz = \left(z , z , \ldots , z\right) : \, z \in \paren{\rr^{d}}^{n} \right\}$. Two important features of this formulation are: (i) the projection onto the set $C$ can be easily computed since
\begin{equation*}
	P_{C}\left(\bz\right) = \left(P_{C_{0}}\left(z_{0}\right) , P_{C_{1}}\left(z_{1}\right) , \ldots , P_{C_{m}}\left(z_{m}\right)\right),
\end{equation*}
where $P_{C_{j}}$, $j = 1 , 2 , \ldots , m$, are given in \eqref{eq:P_Cj}, and (ii) $D$ is a subspace which also has simple projection given by $P_{D}\left(\bz\right) = {\bar \bz}$ where
\begin{equation*}
	{\bar z}_{j} = \frac{1}{m + 1}\sum_{j = 0}^{m} z_{j}.
\end{equation*}
The most natural algorithm for this formulation is the one we began with, namely Cyclic Projections Algorithm \ref{alg:cp}. In the case of just two sets, this is known as {\em Alternating Projections} Algorithm:

    \begin{center}
    		\fcolorbox{black}{Ivory2}{\parbox{15.7cm}{\vspace{-0.1in}
		\begin{alg}[{\bf Alternating Projections - AP}]\label{alg:ap} $~$\\
    			{\bf Initialization.} Choose $\bz^{0} \in \paren{\rr^{d}}^{n(m + 1)}$. \\
    			{\bf General Step ($k = 0 , 1 , \ldots$)}
			    \begin{equation*}
					\bz^{k + 1} \in P_{D}P_{C}\bz^{k}.
			    \end{equation*}
		\end{alg}\vspace{-0.1in}}}
	\end{center}

One can easily verify that Algorithm \ref{alg:ap} is equivalent to the following iteration:
\begin{equation*}
	z_{j}^{k + 1} \in \frac{1}{m + 1}\sum_{j = 0}^{m} P_{C_{j}}\left(z_{j}^{k}\right), \quad j = 0 , 1 , \ldots , m;
\end{equation*}
in other words, Alternating Projections Algorithm on the product space coincides with the Averaged Projections Algorithm \ref{alg:avp} and the Alternating Minimization Algorithm \ref{alg:am}. It is also interesting to note that, in the product space, the Alternating Projections Algorithm is equivalent to the very popular {\em Projected Gradient} Method. To see this, consider the following minimization problem:
\begin{equation}\label{e:duh}
	\ucmin{\frac{1}{2}\dist^{2}\left(\bz , D\right)}{\bz \in C} .
\end{equation}
The objective of this minimization problem is convex and continuously differentiable (since $D$ is a closed and convex set) with a Lipschitz continuous gradient (with constant $1$) given by $\nabla \dist^{2}\left(\bz , D\right) = 2\left(\bz - P_{D}\bz\right)$ (see \eqref{eq:grad sum-of-sq-dist}). The classical Projected Gradient Algorithm applied to this problem follows immediately.

    \begin{center}
    		\fcolorbox{black}{Ivory2}{\parbox{15.7cm}{\vspace{-0.1in}
		\begin{alg}[{\bf Projected Gradient - PG}]\label{alg:pg} $~$\\
    			{\bf Initialization.} Choose $\bz^{0} \in \paren{\rr^{d}}^{n(m + 1)}$. \\
    			{\bf General Step ($k = 0 , 1 , \ldots$)}
			    \begin{eqnarray*}
				\bz^{k + 1} \in P_C\left(\bz^{k} -  \tfrac12\nabla \dist_{D}^{2}(\bz^k)\right) & \iff & \bz^{k + 1} \in P_{C}P_{D}\bz^{k} \\
				& \iff & u^{k + 1} \in \frac{1}{m + 1}\sum_{j = 0}^{m} P_{C_{j}}u^{k}, \,\, \text{where} \, z_{j}^{k + 1} = P_{C_{j}}u^{k + 1}.
			    \end{eqnarray*}
		\end{alg}\vspace{-0.1in}}}
	\end{center}

This is not surprising since the minimization problem \eqref{e:duh} is equivalent to \eqref{SumModel:2}. Indeed, by the definition of the distance function we obtain that
\begin{align*}
	\min_{\bz \in C} \frac{1}{2}\dist^{2}\left(\bz , D\right) & = \min_{\bz \in C}\min_{\bu \in D} \frac{1}{2}\norm{\bz - \bu}^{2} = \min_{\bz \in C}\min_{u \in \paren{\rr^{d}}^{n}} \frac{1}{2}\sum_{j = 0}^{m}\norm{u - z_{j}}^{2} \\
& = \min_{u , \bz} \set{\sum_{j = 0}^{m} \frac{1}{2}\norm{u - z_{j}}^{2}}{z_{j} \in C_{j}, \quad j = 0 , 1 , \ldots , m}.
\end{align*}
It is well-known that the Projected Gradient Algorithm can be accelerated in the convex setting \cite{Nesterov83,BeckTeboulle09}. Although there is no theory that supports acceleration in the nonconvex setting, the recent work \cite{PauBecEldSab18} demonstrates successful empirical results for a class of phase retrieval problems. Following this line, a Fast Projected Gradient Algorithm for problem \eqref{e:duh} reads as follows:

    \begin{center}
    		\fcolorbox{black}{Ivory2}{\parbox{15.7cm}{\vspace{-0.1in}
		\begin{alg}[{\bf Fast Projected Gradient - FPG}]\label{alg:favp} $~$\\
    			{\bf Initialization.} Choose $\bz^{0} , \by^{1}\in \paren{\rr^{d}}^{n(m + 1)}$ and $\alpha_{k} =\frac{k - 1}{k + 2}$ for all $k \in \nn$. \\
    			{\bf General Step ($k = 1 , 2 , \ldots$)}
			    \begin{eqnarray*}
					\bz^{k} & \in & P_{C}\left(\by^{k} -  \tfrac12\nabla \dist^{2}\left(\by^{k} , D\right)\right), \\
                    	\by^{k + 1} & = & \bz^{k} + \alpha_{k}\left(\bz^{k} - \bz^{k - 1}\right) \\
					&\iff&\\
					\bz^{k} &\in & P_{C}P_{D}\by^{k}, \\
                  	\by^{k + 1} & = & \bz^{k} + \alpha_{k}\left(\bz^{k} - \bz^{k - 1}\right). \\
			    \end{eqnarray*}
		\end{alg}\vspace{-0.1in}}}
	\end{center}

There is no theory for the choice of acceleration parameter $\alpha_{k}$, $k \in \nn$, in Algorithm \ref{alg:favp} for nonconvex problems, but our numerical experiments indicate that an investigation into this would be fruitful.

Another algorithmic approach that can fit to the setting of two-set feasibility is the now popular Douglas-Rachford algorithm \cite{DougRach56}. To compensate for the absence of fixed points for inconsistent feasibility, Luke proposed in \cite{Luke05a}, a relaxation of this algorithm, which can be viewed as the usual Douglas-Rachford Algorithm for more general proximal mappings applied to a relaxation of the feasibility problem to a smooth constrained optimization problem
\begin{equation}\label{eq:DRlambda}
\ucmin{ \left\{ \frac{\lambda}{2\left(1 - \lambda\right)}\dist^{2}\left(\bz , D\right) + \iota_{C}\left(\bz\right) \right\}}{\bz \in \paren{\rr^{d}}^{n(m + 1)}}.
\end{equation}

    \begin{center}
    		\fcolorbox{black}{Ivory2}{\parbox{15.7cm}{\vspace{-0.1in}
		\begin{alg}[{\bf Relaxed Douglas-Rachford - DR$\lambda$}]\label{alg:raar} $~$\\
    			{\bf Initialization.} Choose $\bz^{0} \in \paren{\rr^{d}}^{n(m + 1)}$ and $\lambda \in \left[0 , 1\right]$. \\
    			{\bf General Step ($k = 0 , 1 , \ldots$)}
			    \begin{equation}
					\bz^{k + 1} \in \frac{\lambda}{2}\left(R_{D}R_{C}\bz^{k} + \bz^{k}\right) + \left(1 - \lambda\right)P_{C}\bz^{k}.
			    \end{equation}
		\end{alg}\vspace{-0.1in}}}
	\end{center}
	
In \cite{Luke05a} and \cite{Luke08} this algorithm is called the {\em Relaxed Averaged Alternating Reflections} (RAAR) Algorithm. It was shown in \cite{Luke08} that this fixed point mapping is precisely the proximal Douglas-Rachford Algorithm applied to the problem \eqref{eq:DRlambda}, that is,
\begin{equation*}
	\frac12 \paren{R_{1}R_{C} + \Id} = \frac{\lambda}{2}\left(R_{D}R_{C} + \Id\right) + \left(1 - \lambda\right)P_{C},
\end{equation*}
where $R_{1}$ is the \textit{proximal reflector} of the function $f_{\lambda}\left(\bz\right) \equiv \frac{\lambda}{2(1 - \lambda)}\dist^{2}\left(\bz , D\right)$, that is, $R_{1}\left(\bz\right) = 2\prox_{1,f_{\lambda}}\left(\bz\right) - \bz$. The analysis of \cite{LukNguTam17} also applies to the product space formulation of phase retrieval and source localization.

A similar algorithm called Hybrid Projection Reflections was proposed in \cite{BCL2}, but we do not include this in our comparisons because this algorithm, like the original HIO algorithm \cite{Fien82} that inspired it, is not stable.

An interesting alternative to Douglas-Rachford and Alternating Projections Algorithms is the following algorithm, which in a limiting case is simply a convex combination of the two algorithms. For an analysis of this algorithm and a characterization of the set of fixed points see \cite{Thao18}.

    \begin{center}
    		\fcolorbox{black}{Ivory2}{\parbox{15.7cm}{\vspace{-0.1in}
		\begin{alg}[{\bf Douglas-Rachford-Alternating-Projections - DRAP}]\label{alg:drap} $~$\\
    			{\bf Initialization.} Choose $\bz^{0} \in \paren{\rr^{d}}^{n(m + 1)}$ and $\lambda \in \left[0 , 1\right]$. \\
    			{\bf General Step ($k = 0 , 1 , \ldots$)}
			    \begin{equation}
					\bz^{k + 1} \in P_{D}\left(\left(1 + \lambda\right)P_{C}\bz^{k} - \lambda \bz^{k}\right) - \lambda\left(P_{C}\bz^{k} - \bz^{k}\right).
			    \end{equation}
		\end{alg}\vspace{-0.1in}}}
	\end{center}

\section{Numerical Comparisons} \label{s:numerics}
As mentioned in the introduction, there are a lot of algorithms one could choose from for comparisons, and we do not include all here. The algorithms  that simply are not competitive, like {\em Phase Lift} or {\em HPR}, are not included. Other algorithms that are frequently used in the applications literature are either equivalent to one of the algorithms above (HIO \cite{Fien82} and the {\em Difference Map} Method \cite{Elser1} are either Douglas-Rachford or HPR depending on the qualitative constraints; Gerchberg-Saxton \cite{Gerch72} and Error Reduction \cite{Fien82} are the Alternating Projections Algorithm) or are refinements (preconditioning, step-length optimization, etc.) of one of the featured algorithms. Still other methods are too specific to the application to be generalized easily.

Our comparisons report iteration counts. The $^{\ast}$ in the iteration count field of the tables indicates that the algorithm does not converge. The per iteration operation count is comparable between all algorithms except the QNAvP Algorithm (see Algorithm \ref{alg:LBFGS}), where the per iteration cost is approximately five times that of the other algorithms. The failure of the Dynamically Reweighted Averaged Projections (DyRePr) Algorithm \ref{alg:DyRePr} in many of tables below was given an asterisk to note that the fixed points of this algorithm need not be Euclidean best approximation points, so the failure of the algorithm to come within a fixed distance of the true solution may not indicate that the fixed point of this algorithm is bad. Indeed, the DyRePr algorithm converges to a critical point without fail, but as measured by the distance to the true solution, this is not within the specified tolerance, with or without noise. Nonetheless, the fixed points of this algorithm generally {\em looked} 
good to the naked eye, and these points are probably still scientifically informative, even if our rather stringent criteria do not reflect this. It is curious, however, that the objective function underlying this algorithm is the only one motivated by statistical reasoning.

Our results show that the preference for smooth models, while reasonable in the context of the standard (smooth) analytical tools, is not justified by algorithm performance. Moreover, recent advances in the analysis of nonsmooth and nonconvex optimization have placed the Cyclic Projections Algorithm \ref{alg:cp} - even in the inconsistent, nonconvex setting - on as firm theoretical ground as methods for smooth optimization. The Cyclic Relaxed Douglas-Rachford Algorithm \ref{alg:cdr} still remains to be analyzed for the phase retrieval problem, but the path forward is relatively clear.

\subsection{Wavefront Sensing: JWST dataset}
We begin with numerical comparisons on the James Webb Space Telescope wavefront reconstruction data set presented in \cite{Luke02a} (at that time the telescope was known  only as the Next Generation Space Telescope). This is synthetic data, but was carefully constructed to simulate the actual telescope. For our numerical experiments we use a resolution of $128 \times 128$ for each of the images. The data for this problem, $b_{ij} \in \paren{\rr^{2}}$ in problem \eqref{eq:physical model} with $n = 128^{2}$, consists of two out-of-focus images and one in-focus image of a known star ($m = 3$), together with the constraint that the wavefront to be constructed has unit amplitude across the aperture of the telescope (see \cite{Luke02a}). The focus settings and aperture support are accounted for in the mappings $\Pcal_{j}$, $j = 1 , 2 , \ldots , m$, in problem \eqref{eq:physical model}.

All algorithms started from the same randomly chosen initial points and were terminated once the difference between successive iterates falls below $5\times 10^{-5}$, or a maximum iteration count ($6000$) is exceeded. The algorithm is judged to have failed if it does not achieve a distance to the true solution of $10^{-2}$ relative to global phase shift for noiseless data, or a distance of $0.35$ for noisy data.

The Averaged Projections Algorithm \ref{alg:avp}, Dynamically Reweighted Averaged Projections Algorithm \ref{alg:DyRePr} and limited memory BFGS Algorithm \ref{alg:LBFGS} were already compared in \cite{Luke02a}.

An interesting feature of this data set that none of the other experiments share is that the Cyclic Projections (CP) and the Cyclic Relaxed Douglas-Rachford (CDR$\lambda$) Algorithms do not converge in the sense we might expect. Instead, what the iterates converge to is a fixed amplitude and fixed {\em relative phase}, but with a constant {\em global} phase shift from one iterate to the next: each pixel is like a marble spinning forever around and around, in synchrony with all the other pixels. What we observe is linear convergence of the iterates rotated by this (limiting) global phase step. Also, the behavior of CDR$\lambda$ is sensitive to the choice of $\lambda$ in the same way DR$\lambda$ is: if the parameter is too large, the algorithm will not converge in any sense. This behavior has yet to be explained.

\begin{table}[ht]
\begin{center}{\tiny
\begin{tabular}{||c||c|c|c|c||c|c|c|c|}\hline
& \multicolumn{4}{c||}{$~$}&\multicolumn{4}{c|}{$~$}\\
& \multicolumn{4}{c||}{{\large JWST no noise}}&\multicolumn{4}{c|}{{\large JWST noise}} \\
& \multicolumn{4}{c||}{iteration count}&\multicolumn{4}{c|}{iteration count}
\\
$~$&failure & median  & high& low & failure & median & high & low \\\hhline{||=||=|=|=|=||=|=|=|=|}
Wirtinger & $100$ & $86$ &  $87$ & $85$
& $100$ & $86$ &  $87$ & $85$
\\\hhline{||-||-|-|-|-||-|-|-|-|}
ADMM$_1$ ($\eta=3$) & $0$ & $70.5$ &  $115$ & $47$
& $100$ & $*$ &  $*$ & $*$
\\\hhline{||-||-|-|-|-||-|-|-|-|}
ADMM$_2$ ($\rho_j=.5$)& $4$ & $1272$ &  $6000$ & $696$
& $0$ & $917.5$ &  $5201$ & $565$
\\\hhline{||-||-|-|-|-||-|-|-|-|}
AP/AvP/PG & $0$ & $247$ &  $1700$ & $95$
& $0$ & $138$ &  $508$ & $77$
\\\hhline{||-||-|-|-|-||-|-|-|-|}
DyRePr & $100^*$ & $221$ &  $2039$ & $121$
& $100^*$ & $151.5$ &  $586$ & $88$
\\\hhline{||-||-|-|-|-||-|-|-|-|}
FPG & $0$ & $191$ &  $789$ & $87$
& $0$ & $145.5$ &  $589$ & $80$
\\\hhline{||-||-|-|-|-||-|-|-|-|}
QNAvP & $3$ & $95$ &  $240$ & $47$
& $3$ & $67.5$ &  $206$ & $34$
\\\hhline{||-||-|-|-|-||-|-|-|-|}
DR & $68$ & $704.5$ &  $821$ & $631$
& $100$ & $*$ &  $*$ & $*$
\\\hhline{||-||-|-|-|-||-|-|-|-|}
DR$\lambda$ ($\lambda=0.85/0.55$) & $0$ & $64.5$ &  $145$ & $53$
& $0$ & $109$ &  $515$ & $62$
\\\hhline{||-||-|-|-|-||-|-|-|-|}
DRAP ($\lambda=0.55/0.25$) & $0$ & $72.5$ &  $319$ & $45$
& $0$ & $97$ &  $362$ & $57$
\\\hhline{||-||-|-|-|-||-|-|-|-|}
CP$^*$  & $0$ & $27$ &  $201$ & $17$
& $0$ & $21$ &  $48$ & $14$
\\\hhline{||-||-|-|-|-||-|-|-|-|}
CDR  & $0$ & $9$ &  $18$ & $8$
& $100$ & $*$ &  $*$ & $*$
\\\hhline{||-||-|-|-|-||-|-|-|-|}
CDR$\lambda^*$ ($\lambda = 0.25$)  & $0$ & $25$ &  $423$ & $15$
& $0$ & $21$ &  $81$ & $14$ \\
\hline
\end{tabular}}
\caption{\label{table JWST} {\small JWST wavefront reconstruction problem, noiseless and noisy, $100$ random initializations. The asterisk on $CP$ and $CDR\lambda$ are to note that these algorithms with noisy JWST data {\em do not converge} to fixed points, but rather to a fixed amplitude with a constant global phase shift from one iterate to the next. The $^{\ast}$ for DyRePr indicates that the algorithm fixed point still {\em looks} good, even if it does not lie within the specified distance to the true solution.}}
\end{center}
\end{table}

\subsection{Coded Diffraction: CDP Dataset}
A phase problem for demonstrating the Wirtinger Flow Algorithm was presented in \cite{candes2014phase}.
The experiments are of synthetic $1$- and $2$-dimensional signals. We compare the algorithms surveyed here on this problem instance for the same data made available in \cite{candes2014phase}. There are several features of this problem that have attracted attention in the applied mathematics and statistics communities in recent years. Regarding the physical experiment that this data set mimics, it is imagined that one has a {\em phase mask} at the object plane that one can randomly generate. The data consists of $10$ observations of a true signal, each observation $b_{ij} \in \rr^{2}$ in problem \eqref{eq:physical model} with $n = 128^{2}$ and $m = 10$ made with a different, randomly generated, phase mask - $\Pcal_{j}$, $j = 1 , 2 , \ldots , m$, in problem \eqref{eq:physical model}. There is no qualitative constraint set $C_{0}$. To avoid getting stuck in a bad local minimum, the Wirtinger Flow approach involves a {\em warm start} procedure which is a power series iteration ($50$ in our experiments) that is 
meant to land one in the neighborhood of a global minimum. The warm start strategy proposed in \cite{candes2014phase} did not have any significant beneficial effect on avoidance of bad local minimums - in fact, it made things {\em worse} in several cases. What distinguishes the best performers from the worst for these experiments is not only their computational efficiency (by a factor of $100$) but also the fact that they never failed to find the globally optimal solution with or without a warm start.

This experiment also includes a one-dimensional phase retrieval benchmark. It is held that $1$-dimensional phase problems of this type are very different than $2$-dimensional instances. This a reasonable claim due to the early theoretical work showing that the $1$-dimensional phase retrieval problem suffers from nonuniqueness, while $2$-dimensional phase retrieval problems are almost always uniquely solvable, when they are solvable at all \cite{Bruck79}. The distinction, however, is spurious when the phase problem is {\em overdetermined}, as it is here. The numerical tests bear this out: there does not appear to be any qualitative computational difference between $1$- and $2$-dimensional phase retrieval problems when both are overdetermined. The exceptions to this are the algorithms based on the (unrelaxed) Douglas-Rachford mapping (DR and CDR in Table \ref{table cdp}), which appear to be unstable for the $2$-dimensional signals where for the 1-dimensional signals they are stable. This could indicate larger 
basins of attraction for $1$-dimensional compared to $2$-dimensional phase retrieval, but otherwise does not point to any more interesting qualitative difference that we can imagine.

All algorithms started from the same randomly chosen initial points and  were terminated once the difference between successive iterates falls below $1e-10$ ($1$-D) or $1e-8$ ($2$-D) respectively, or a maximum iteration count ($6000$) is exceeded. The algorithm is judged to have failed if it does not achieve a distance of $1e-9$ to the true solution for $1$-D signals or $1e-7$ for $2$-D signals, up to global phase shift, before the iterate difference tolerance is reached. Again, the DyRePr Algorithm is reliable, but has fixed points that are not close enough to the true solution for these comparisons.

\begin{table}[ht]
\begin{center}{\tiny
\begin{tabular}{||c||c|c|c|c||c|c|c|c|}\hline
& \multicolumn{4}{c||}{$~$}&\multicolumn{4}{c|}{$~$}\\
& \multicolumn{4}{c||}{{\large 1D cold/warm}}&\multicolumn{4}{c|}{{\large 2D cold/warm}} \\
& \multicolumn{4}{c||}{iteration count}&\multicolumn{4}{c|}{iteration count}
\\
$~$&failure & median  & high& low & failure & median & high & low \\\hhline{||=||=|=|=|=||=|=|=|=|}
Wirtinger & $15$/$9$ & $322$/$297$ &  $6000$/$6000$ & $290$/$269$
& $0$/$12$ & $412.5$/$266.5$ &  $999$/$404$ & $332$/$6$
\\\hhline{||-||-|-|-|-||-|-|-|-|}
ADMM$_2$ ($\rho_j=.8$)& $0$/$0$ & $1040$/$504.5$ &  $1570$/$1215$ & $824$/$435$
& $0$/$0$ & $1034$/$1050$ &  $1831$/$1814$ & $814$/$799$
\\\hhline{||-||-|-|-|-||-|-|-|-|}
AP/AvP/PG & $0$/$0$ & $131.5$/$99$ &  $217$/$121$ & $107$/$89$
& $0$/$0$ & $229$/$122$ &  $334$/$227$ & $179$/$80$
\\\hhline{||-||-|-|-|-||-|-|-|-|}
DyRePr & $0$/$0$ & $54$/$53$ &  $60$/$59$ & $49$/$47$
& $100^*$/$100^*$ & $202.5$/$91.5$ &  $366$/$390$ & $151$/$47$
\\\hhline{||-||-|-|-|-||-|-|-|-|}
FPG & $0$/$0$ & $86$/$84$ &  $102$/$98$ & $78$/$76$
& $0$/$0$ & $3178$/$1375.5$ &  $7596$/$5866$ & $714$/$156$
\\\hhline{||-||-|-|-|-||-|-|-|-|}
QNAvP & $0$/$0$ & $25$/$24$ &  $28$/$28$ & $21$/$20$
& $3$/$12$ & $247$/$52.5$ &  $616$/$230$ & $150$/$28$
\\\hhline{||-||-|-|-|-||-|-|-|-|}
DR & $0$/$0$ & $*$/$*$ &  $*$/$*$ & $131$/$136$
& $100$/$100$ & $*$/$*$ &  $*$/$*$ & $*$/$*$
\\\hhline{||-||-|-|-|-||-|-|-|-|}
DR$\lambda$ ($\lambda=0.75$) & $0$/$0$ & $97$/$91$ &  $105$/$97$ & $92$/$87$
& $0$/$0$ & $88$/$79$ &  $97$/$87$ & $85$/$73$
\\\hhline{||-||-|-|-|-||-|-|-|-|}
DRAP ($\lambda=0.55$) & $0$/$0$ & $68.5$/$60$ &  $86$/$67$ & $63$/$58$
& $0$/$0$ & $57$/$58$ &  $93$/$72$ & $68$/$49$
\\\hhline{||-||-|-|-|-||-|-|-|-|}
CP  & $0$/$0$ & $16$/$13$ &  $20$/$15$ & $13$/$12$
& $0$/$0$ & $19$/$15$ &  $26$/$19$ & $17$/$12$
\\\hhline{||-||-|-|-|-||-|-|-|-|}
CDR  & $0$/$2$ & $39$/$38.5$ &  $6000$/$6000$ & $20$/$18$
& $100$/$100$ & $*$/$*$ &  $*$/$*$ & $*$/$*$
\\\hhline{||-||-|-|-|-||-|-|-|-|}
CDR$\lambda$ ($\lambda = 0.33$)  & $0$/$0$ & $6$/$6$ &  $6$/$6$ & $5$/$5$
& $0$/$0$ & $11$/$9$ &  $14$/$12$ & $11$/$8$ \\
\hline
\end{tabular}}
\caption{\label{table cdp} {\small
\noindent $1$- and $2$-dimensional phase retrieval, problem CDP \cite{candes2014phase}, cold/warm start, $100$ random initializations. The $^{\ast}$ for DyRePr indicates that the algorithm fixed point still yields a point close to the true solution, though it does not lie within the specified distance for this comparison.}}
\end{center}
\end{table}

\subsection{Sparse Phase retrieval}
The first to examine sparsity constrained phase retrieval was Marchesini in \cite{Marchesini08} where he recognized the  implicit use of sparsity and hard thresholding in phase retrieval in the {\em Charge Flipping} Algorithm of Oszlanyi and S\"ut\'o \cite{Oszlanyi04,Oszlanyi08}. Since then, sparse phase retrieval has received intense study. More than one quarter of the results returned from a web query of the term ``phase retrieval'' are about sparse phase retrieval. Our contribution to this wave includes the recent studies \cite{Thao18b, PauBecEldSab18,BSTV2018}. All of the algorithms surveyed in the previous section have analogues as more general iterated prox-algorithms (the projection onto a set is just the prox of the indicator function) and so the projections can be replaced quite easily with hard and soft-thresholders which are, respectively, the prox mappings of the counting function ($\norm{x}_{0} = \sum_{j = 1}^{n} \norm{x_{j}}_{0}$) and the $\ell_{1}$-norm.

The sparsity example constructed for our experiment is modeled after \cite{Marchesini18}. We generate a fixed number of Gaussian densities of different height and width randomly placed on a $128 \times 128$ discretization of the plane. The data, $b_{ij} \in \rr^{128^{2}}$ in problem \eqref{eq:physical model} with $n = 128^{2}$ and $m = 1$, is the amplitude of the discrete Fourier transform of the collection of Gaussian distributions. Instead of promoting sparsity using the counting function or the $\ell_{1}$-norm, we specify an estimate of the maximum number, $s$, of nonzero entries (which for our experiments was an overestimate of the true sparsity by $20\%$) and project onto the set
\begin{equation*}
    S_{s} \equiv \left\{ z \in \paren{\rr^{2}}^{128^{2}} ~|~ \norm{z}_{0} \leq s \right\}.
\end{equation*}
To this we add the constraint that the object is real and nonnegative, i.e., belongs to
\begin{equation} \label{eq:C+}
   C_{+} \equiv \left\{ z = \left(z_{1} , z_{2} , \ldots , z_{128^{2}}\right) \in \paren{\rr^2}^{128^{2}}~|~ z_{j} = \left(x_{j} , 0\right) \in \rr^{2}, ~ x_{j} \in \left[0 , \infty\right) ~ (1 \leq j \leq 128^{2}) \right\}.
\end{equation}
The qualitative constraint is then
\begin{equation*}
   C_{0} \equiv C_{+} \cap S_{s}.
\end{equation*}
It is easy to show that $C_{0}$ is a cone, so this constraint fits the feasibility format of the cone and
sphere problem. Moreover, it is easy to show that $P_{S_{s}}P_{C_{+}} = P_{S_{s} \cap C_{+}}$. As such,
the analysis for the other instances of this problem demonstrated above also apply here. The constraint set $S_{s}$ was employed in \cite{Thao18b}, but to our knowledge, this is the only place where sparsity constrained phase retrieval is modeled in this way, and combination with a real nonnegativity constraint as we have done here appears to be novel.  It is possible to incorporate real-nonnegativity constraints into the hard and soft-thresholding operators, but since these operators showed no advantage over our formulation here, we do not include this comparison.

For this experiment we only compare the most successful of the algorithms from the previous experiments. The goal here is to compare robustness of the best algorithms against local minimums. The algorithm was judged to have succeeded if it correctly identified the support of the Gaussian dots (modulo reflections and translations -- phase retrieval is not unique) up to an error of $5 \times 10^{-4}$.

\begin{figure}[!htp]
  \centering
     \includegraphics[width=11cm, height=5cm]{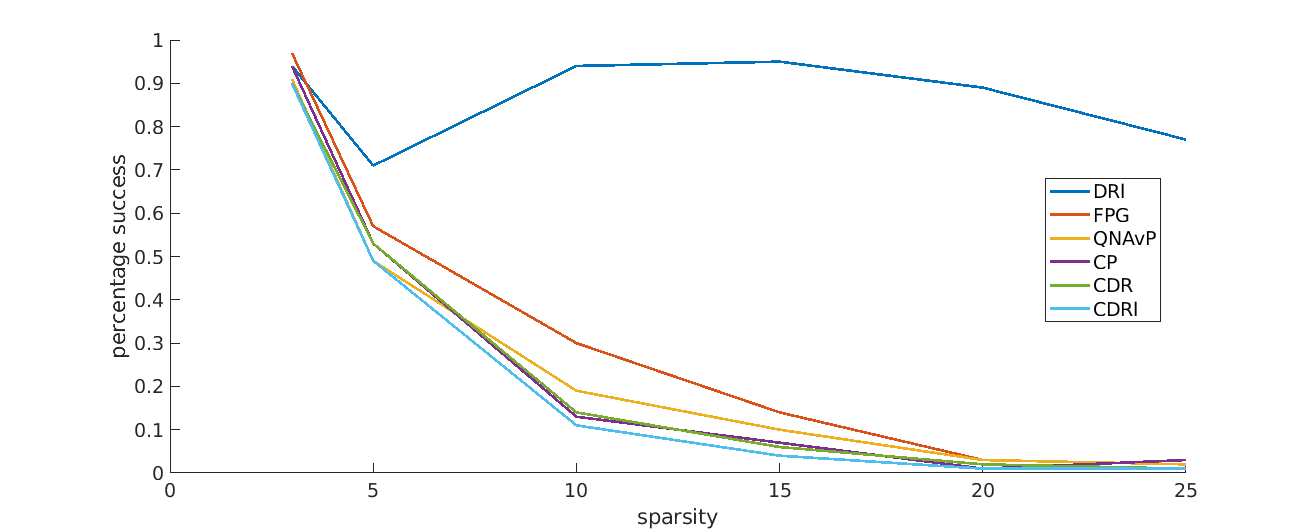}
     \caption{Leading algorithms applied to the feasibility formulation of sparse phase retrieval with nonnegative real sparse objects. Shown are the percentage of successful exact recoveries of the support of the sparse object (up to shifts and reflections) over $100$ random trials for objects with $3, 5, 10, 15, 20$ and $25$ Gaussian ``dots''.}
\end{figure}

What we observe is that there is not a great deal of variability between algorithms with respect to the success rates, with the exception of the Relaxed Douglas-Rachford (DR$\lambda$) Algorithm  \ref{alg:raar}. This is a two-set feasibility problem, so the advantage of the cyclic algorithms CDR and CDR$\lambda$ is not clear. What is not reflected in these graphs is the average number of iterations each algorithm required to reach convergence (they all converged). The relative performance followed the same pattern as the other experiments above, with efficiencies differing by up to two orders of magnitude. The algorithms not included in the comparison did not successfully recover the support and so were not considered.

\subsection{Source Localization}
The sensors are located on a $100 \times 100$ unit grid, and for the noisy experiments the placement of the sensors has an error of about 3 units (SNR $30$). The data consists of $m = 3$ or $10$ measurements $b_{1j} \in \rr$, $j = 1 , 2 , \ldots , m$, with corresponding shift operators $\Pcal_{j}$, $j = 1 , 2 , \ldots , m$, defined by the randomly determined locations of the sensors (see problem \eqref{Ranges}). All algorithms are started from the same randomly chosen initial points. All algorithms are terminated once the difference between successive iterates falls below $1e-11$ or a maximum iteration count ($10000$) is exceeded. The algorithm is judged to have failed if it does not achieve a distance of $1e-1$ to the true source location for noiseless data and $3$ for noisy signals before the iterate difference tolerance is reached. The abstract the source localization problem has the same qualitative geometry as phase retrieval, and the numerical results do not indicate any differences in the relative 
quantitative performance of the algorithms. It is remarkable that the Douglas-Rachford Algorithm converges sometimes to a {\em fixed point} with $3$ sensors and noise. We know that Douglas-Rachford does not possess a fixed point for inconsistent feasibility problems, so this indicates that even with noise, the three sensor case can still have nonempty intersection. Even when the intersection is empty, however, note that the iterates of the Douglas-Rachford Algorithm, while they do not converge, they never wander very far from an acceptable solution. The Wirtinger Flow Algorithm \ref{alg:Wirt} does not make sense in this context so has been dropped from this comparison.

\begin{table}[ht]
\begin{center}{\tiny
\begin{tabular}{||c||c|c|c|c||c|c|c|c|}\hline
& \multicolumn{4}{c||}{$~$}&\multicolumn{4}{c|}{$~$}\\
& \multicolumn{4}{c||}{{\large 3 sensors no noise/noise}}&\multicolumn{4}{c|}{{\large 10 sensors no noise/noise}} \\
& \multicolumn{4}{c||}{iteration counts}&\multicolumn{4}{c|}{{iteration counts}}
\\
$~$&failure & median  & high& low & failure & median & high & low \\\hhline{||=||=|=|=|=||=|=|=|=|}
ADMM$_2$ ($\rho_j=.3$)& $13$/$100$ & $3724.5$/$3723$ &  $10000$/$10000$ & $274$/$302$
& $0$/$90$ & $550.5$/$569.5$ &  $2141$/$10000$ & $284$/$273$
\\\hhline{||-||-|-|-|-||-|-|-|-|}
AP/AvP/PG & $0$/$0$ & $1001.5$/$969$ &  $10000$/$10000$ & $68$/$75$
& $0$/$0$ & $141.5$/$146.5$ &  $570$/$722$ & $71$/$68$
\\\hhline{||-||-|-|-|-||-|-|-|-|}
DyRePr & $100^*$/$0$ & $7$/$6.5$ &  $23$/$23$ & $5$/$5$
& $100^*$/$0$ & $7$/$7$ &  $16$/$16$ & $5$/$5$
\\\hhline{||-||-|-|-|-||-|-|-|-|}
FPG & $0$/$0$ & $847.5$/$961$ &  $10000$/$10000$ & $97$/$103$
& $0$/$0$ & $168.5$/$182.5$ &  $518$/$764$ & $93$/$100$
\\\hhline{||-||-|-|-|-||-|-|-|-|}
QNAvP & $1$/$0$ & $47.5$/$38$ &  $995$/$956$ & $14$/$5$
& $0$/$0$ & $17$/$33$ &  $889$/$46$ & $12$/$13$
\\\hhline{||-||-|-|-|-||-|-|-|-|}
DR & $0$/$0$ & $2248$/$2320$ &  $10000$/$10000$ & $140$/$159$
& $0$/$0$ & $193.5$/$10000$ &  $1197$/$10000$ & $143$/$10000$
\\\hhline{||-||-|-|-|-||-|-|-|-|}
DR$\lambda$ ($\lambda=0.85/0.5$) & $0$/$0$ & $215.5$/$2825$ &  $10000$/$10000$ & $90$/$241$
& $0$/$0$ & $132.5$/$146.5$ &  $209$/$722$ & $102$/$68$
\\\hhline{||-||-|-|-|-||-|-|-|-|}
DRAP ($\lambda=0.55$) & $0$/$0$ & $761$/$498.5$ &  $10000$/$10000$ & $47$/$45$
& $0$/$0$ & $104$/$107.5$ &  $431$/$547$ & $49$/$47$
\\\hhline{||-||-|-|-|-||-|-|-|-|}
CP  & $0$/$0$ & $183.5$/$169.5$ &  $10000$/$7867$ & $9$/$11$
& $0$/$0$ & $10.5$/$11$ &  $55$/$61$ & $5$/$5$
\\\hhline{||-||-|-|-|-||-|-|-|-|}
CDR  & $0$/$0$ & $65$/$68$ &  $10000$/$2871$ & $8$/$8$
& $0$/$0$ & $9$/$9$ &  $54$/$51$ & $4$/$4$
\\\hhline{||-||-|-|-|-||-|-|-|-|}
CDR$\lambda$ ($\lambda = 0.33$)  & $0$/$0$ & $150$/$149$ &  $10000$/$6591$ & $8$/$10$
& $0$/$0$ & $9$/$9$ &  $43$/$52$ & $4$/$4$ \\
\hline
\end{tabular}}
\caption{\label{table sl} {\small
\noindent $100$ random instances of $3$ and $10$ sensor source location problems, with and without noise.  The $^{\ast}$ for DyRePr indicates that the algorithm fixed point is still reasonably close to the truth, even if it does not lie within the specified distance tolerance.}}
\end{center}
\end{table}

\subsection{Experimental Data}
\begin{figure}[!htp]
  \centering
    (a) \includegraphics[width=7cm, height=5cm]{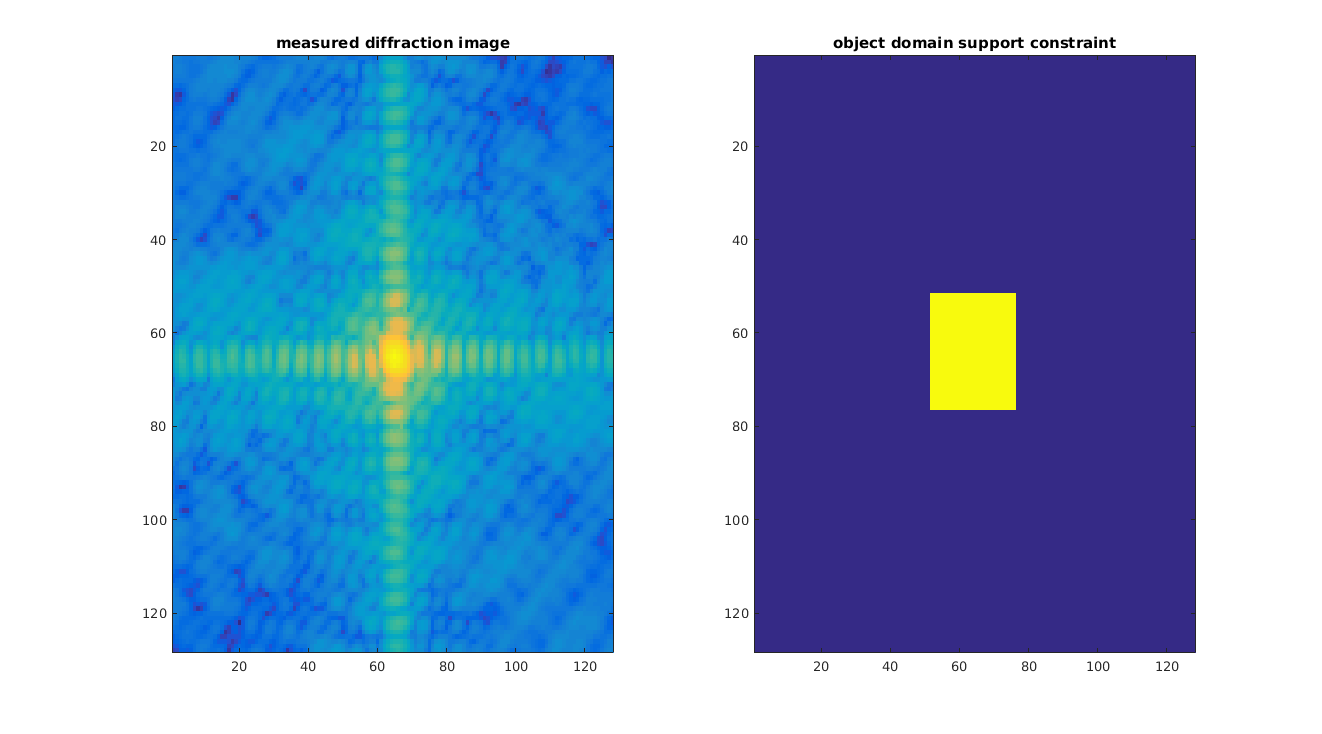}
    \includegraphics[width=7cm, height=5cm]{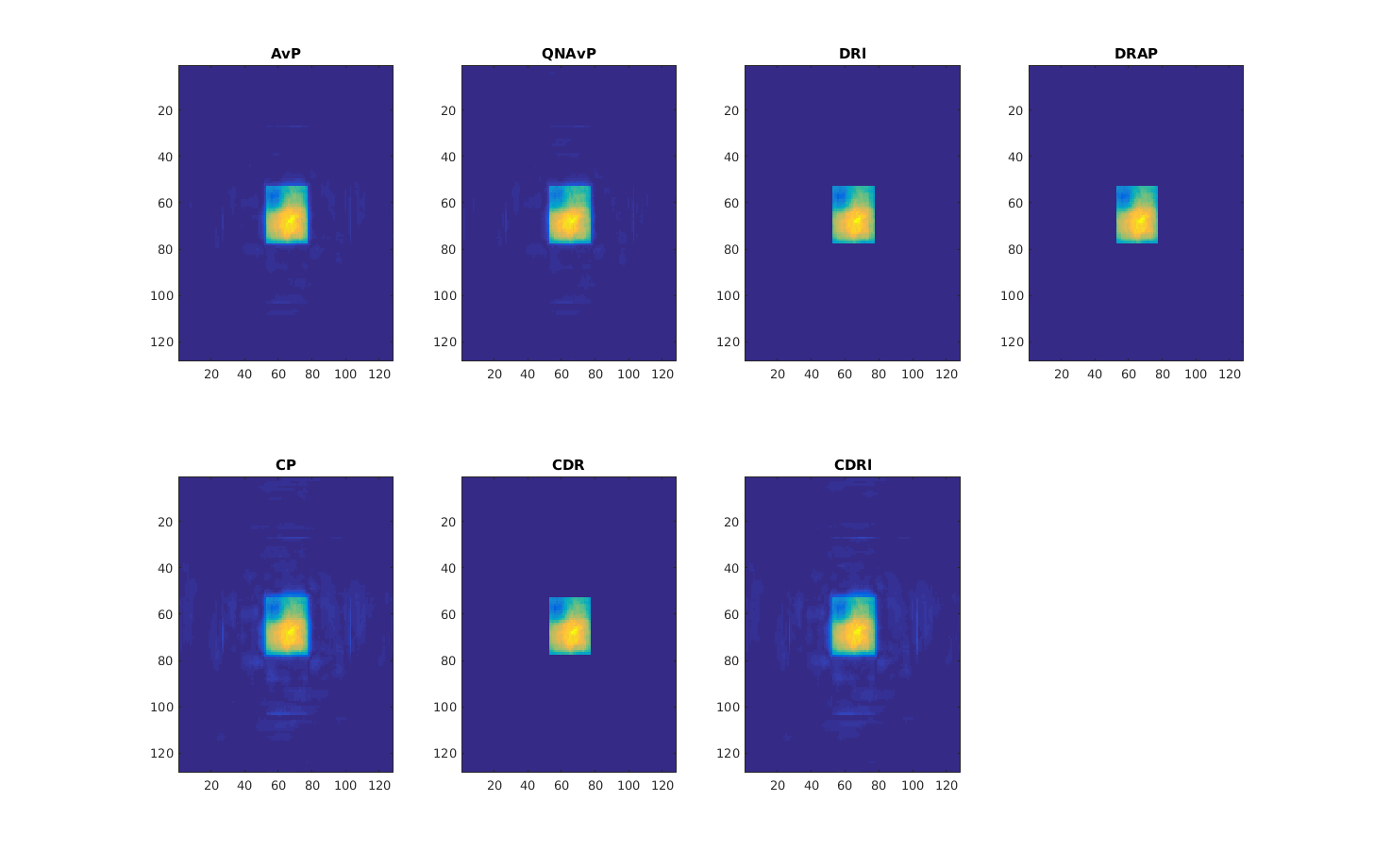}\\
    (b) \includegraphics[width=7cm, height=3cm]{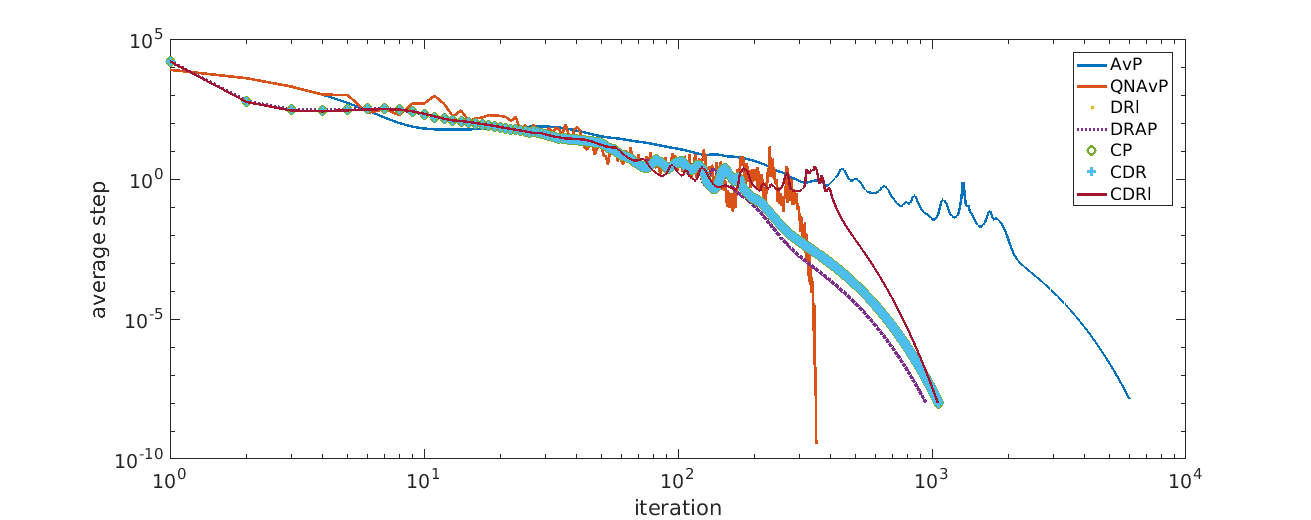}
     \includegraphics[width=7cm, height=3cm]{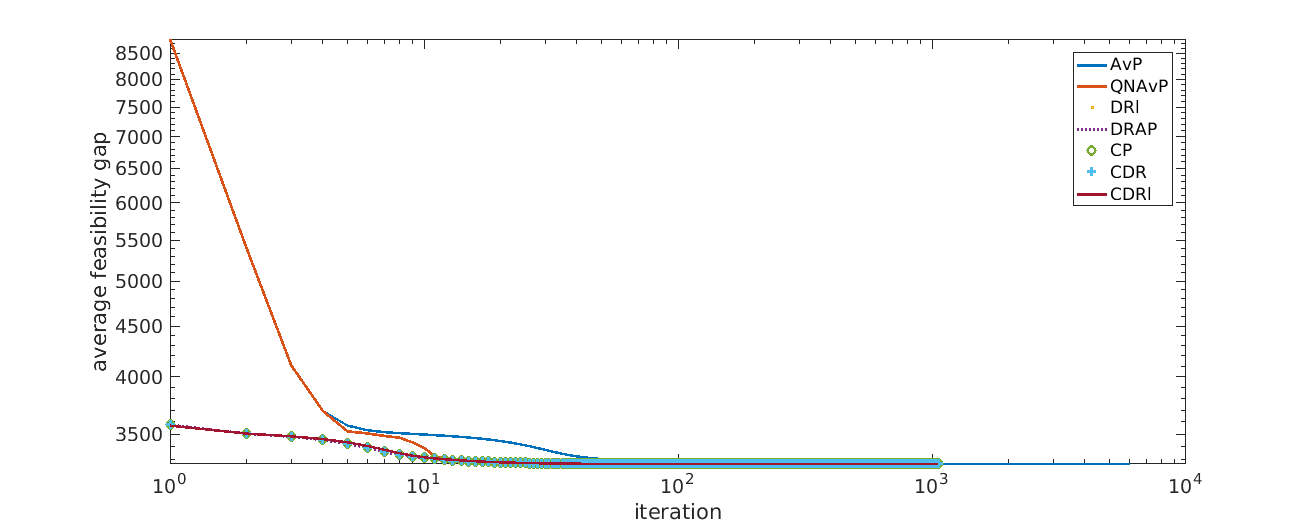}
     \caption{\label{f:tasse} (a) Observation and recovery of an optical diffraction experiment. (b)  Step-size and gap size between constraint sets versus iteraton for several algorithms.}
\end{figure}

We compare the performance of the most successful algorithms above on experimental data. The data set presented here ($b_{ij}$, $j = 1 , 2 , \ldots , m$, in problem \eqref{eq:physical model} with $m = 1$ and $n = 128^{2}$) is an optical diffraction image produced by undergraduates at the X-Ray Physics Institute at the University of G\"ottingen shown in Figure \ref{f:tasse}(a) left. This is a difficult data set because it is very noisy, and the physical parameters of the image (magnification factor, Fresnel number, etc.) were unknown to us. We optimistically assumed a perfect imaging system so that the imaging model is simply an unmodified Fourier transform. We were told that object was a real, nonnegative object, supported on some patch in the object plane, that is, the qualitative constraint $C_{0}$ is of the form \eqref{eq:C+}. In such an experiment, one only has the successive iterates and a feasibility gap to observe in order to conclude that the algorithm is converging at least to a local best 
approximation point. Given the noise in the image, it is not even clear that one desires a local best approximation point with the smallest feasibility gap since this will also mean that the noise has been recovered. For the numerical test reported here we applied a low-pass filter to the data since almost all of the recoverable information about the object was contained in the low-frequency elements. This also had the numerically beneficial effect of rendering the problem {\em more inconsistent}. Thanks to the analysis in \cite{LukNguTam17} we now understand why this can be helpful. The original object was a coffee cup which the generous reader can see if he tilts his head to the left and squints real hard Figure \ref{f:tasse}(a) right.

Noisy experiments like this one also bring to the foreground the issue of {\em regularization}. In \cite{Luke12} approximate projection algorithms for regularized sets were analyzed for cyclic projections.  The ideas of that work can be extended to the other algorithms considered here, namely, that one can place balls (either Euclidean or Kullback-Leibler, as appropriate) around the measured data and project onto these ``inflated'' sets.  Since the projections onto such sets is generally much more complicated to compute than the original sets, one can replace the exact projection with the projection onto the original, unregularized data set, treating the latter projection as an approximation to the former.  This has additional advantages of {\em extrapolation} which yields {\em finite termination} at a feasible points, in this case.

\begin{figure}[!htp]
     (a) \includegraphics[width=15cm, height=5cm]{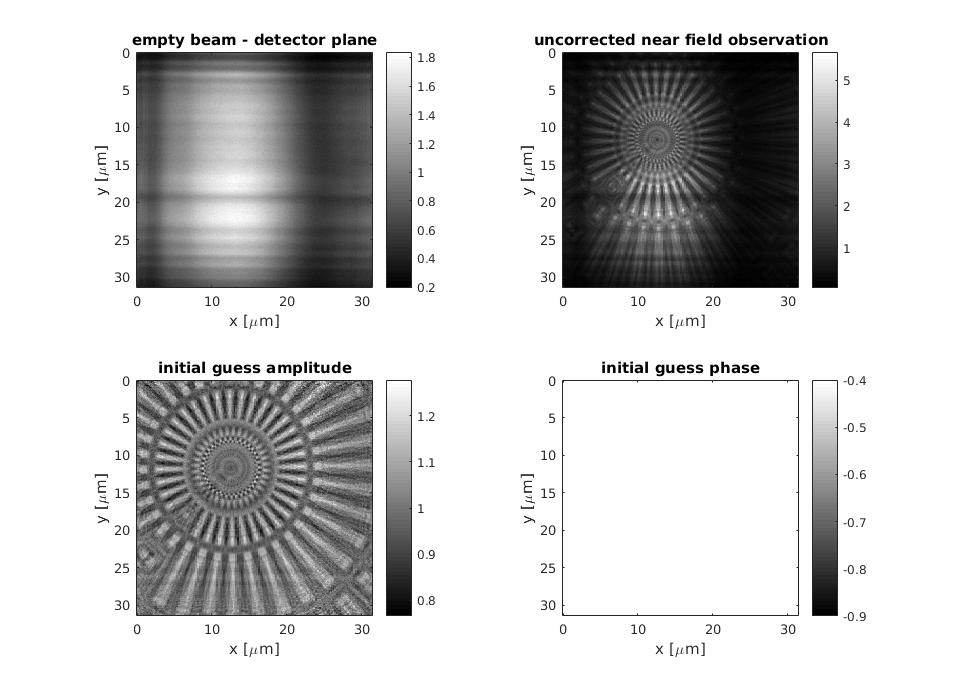}\\
     (b) \includegraphics[width=7cm, height=2cm]{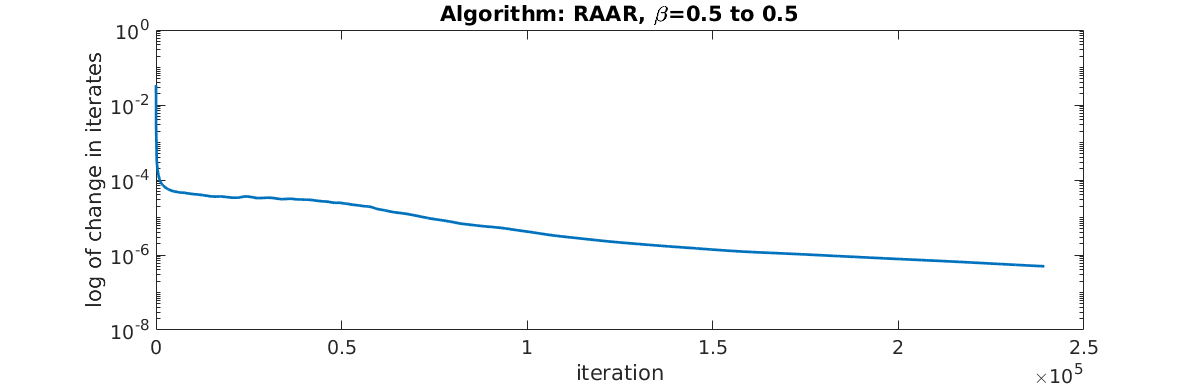}
     \includegraphics[width=7cm, height=2cm]{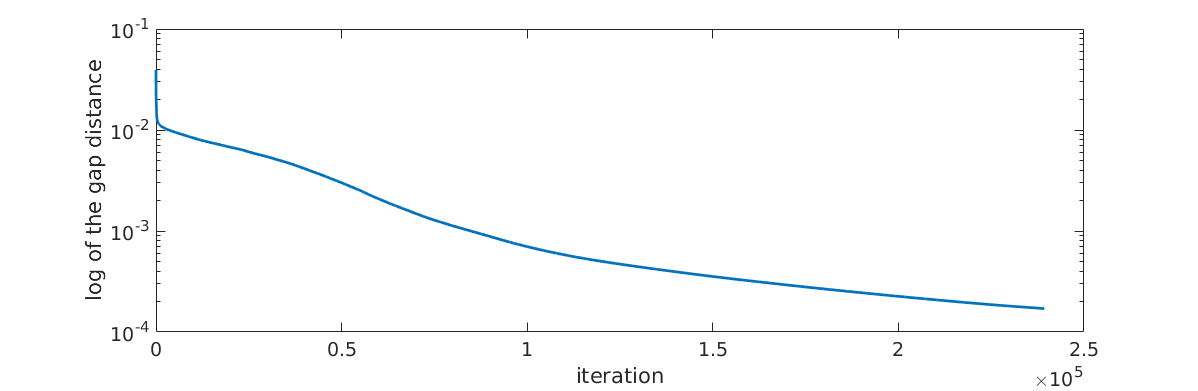}\\
     (c) \includegraphics[width=17cm, height=5cm]{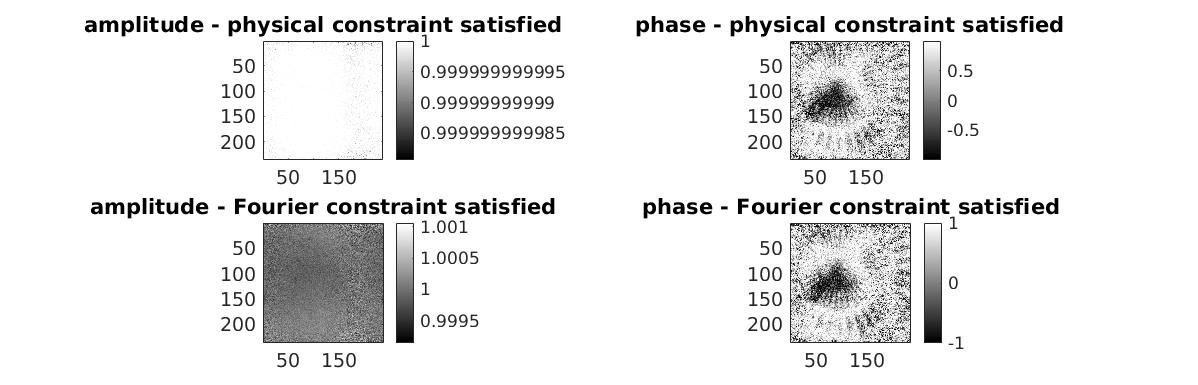}
\caption{\label{fig:Siemens200K} Unregularized reconstruction of near field holography experiment with empty beam correction using the DR$\lambda$ algorithm with $\lambda=0.5$.  The algorithm is stable, and converges as predicted by the theory, however the noise is reconstructed without regularization of the data set.}
\end{figure}

\begin{figure}[!htp]
     (a) \includegraphics[width=15cm, height=5cm]{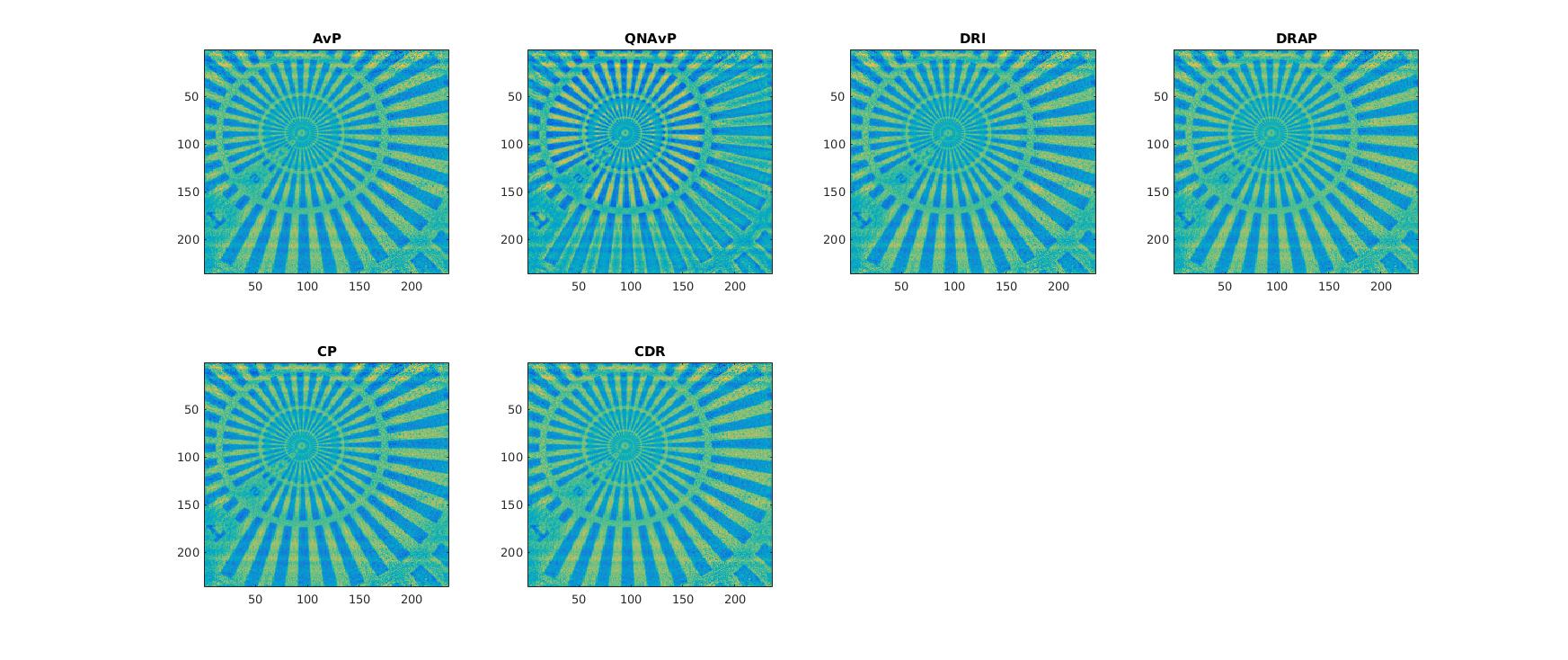}\\
     (b) \includegraphics[width=8cm, height=5cm]{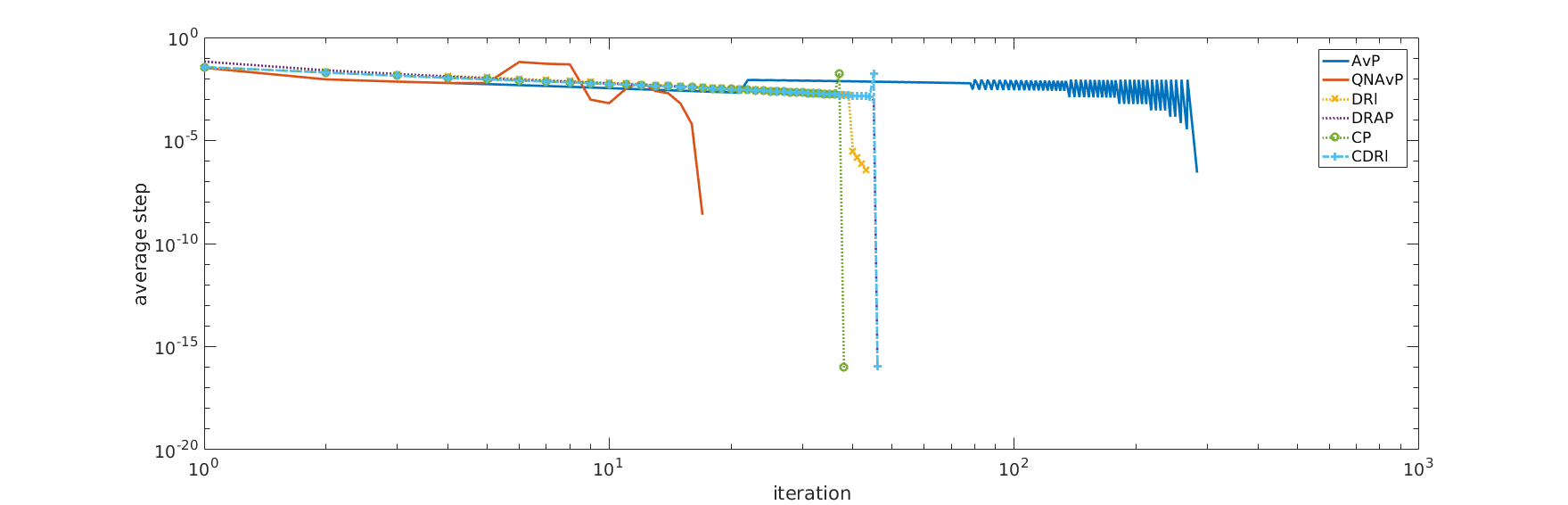}
     \includegraphics[width=8cm, height=5cm]{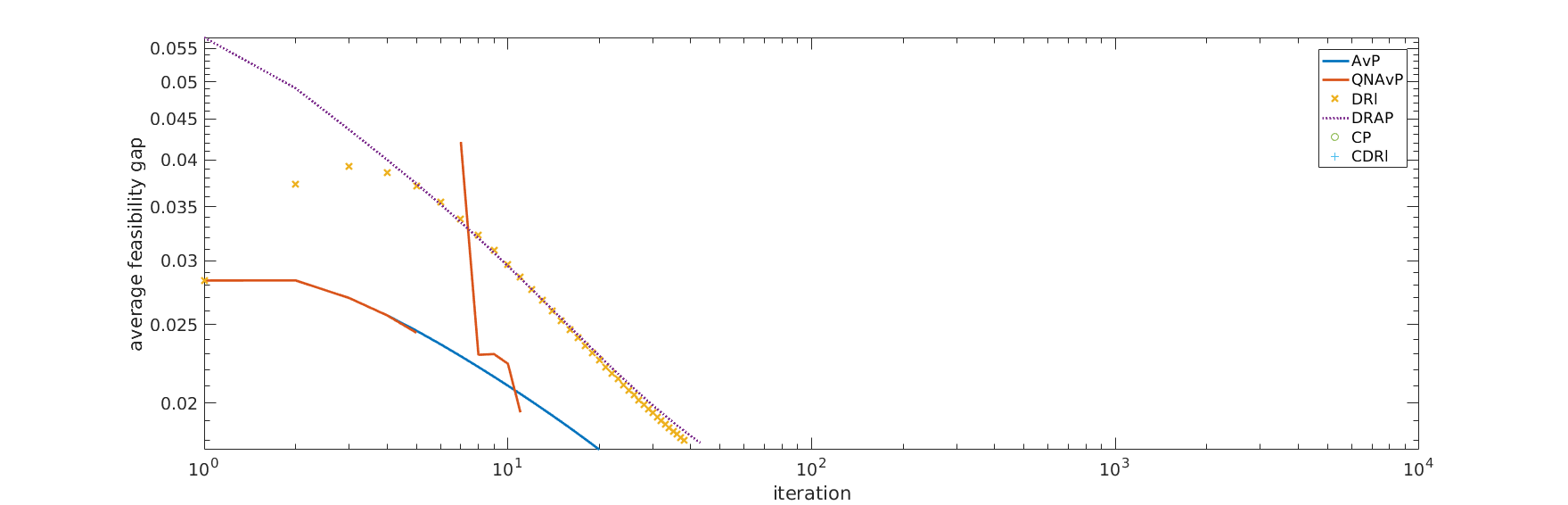}
\caption{\label{fig:Siemens_reg} (a) Reconstruction of regularized near field holography experiment with
empty beam correction for the same data shown in Figure \ref{fig:Siemens200K}(a). (b) Step-size and gap between the 5 constraint sets versus iteration.   }
\end{figure}

The payoff for early termination is demonstrated in Figures \ref{fig:Siemens200K} and \ref{fig:Siemens_reg} with a near field holography experiment provided to us by Tim Salditt's laboratory at the Institute for X-Ray Physics at the University of G\"ottingen. Here the structured illumination shown in Figure \ref{fig:Siemens200K}(a) left is modeled by $\Pcal_{j}$, $j = 1 , 2 , \ldots , m$, in problem \eqref{eq:physical model} with $m = 1$. The image -- $b_{ij}$, $j = 1 , 2 , \ldots , m$, in problem \eqref{eq:physical model} with $m = 1$ --  shown in Figure \ref{fig:Siemens200K}(a) right is in the near field, so the mapping $\Fcal$ in problem \eqref{eq:physical model} is the {\em Fresnel transform} \cite{Hagemann14}. The qualitative constraint is that the field in the object domain has amplitude $1$ at each pixel, that is, the object is a {\em pure phase object}. A support constraint is not applied. Without regularization/early termination the noise is recovered as shown in Figure \ref{fig:Siemens200K}(c). The 
data is regularized by accepting points within a fixed pointwise distance to the measured data at each pixel with respect to the Kullback-Leibler divergence. Rather than projecting onto these Kullback-Leibler balls, all algorithms compute the unregularized projection and move to the point given by the unregularized algorithm. The stopping rule, however, is with respect to the achieved feasibility gap between the regularized data sets and the qualitative constraint. This is effectively an early stopping rule for the unregularized algorithms. The result for different algorithms is shown in Figure \ref{fig:Siemens_reg}.

%

\end{document}